\documentclass[12pt]{article}
\usepackage{amsmath, amssymb,amsfonts}
\usepackage{graphicx}
\usepackage{color}
\providecommand{\U}[1]{\protect\rule{.1in}{.1in}}
\newtheorem{teo}{Theorem}
\newtheorem{coro}{Corollary}
\newtheorem{lemma}{Lemma}

\newtheorem{definition}{Definition}
\newtheorem{remark}{Remark}

\begin{document}

\title{Euler equations with  non-homogeneous Navier slip boundary condition}
\author{N.V. Chemetov* and S.N. Antontsev** \\
{\small }\\
{\small  * DCM-FFCLRP / University of Sao Paulo, Brazil }\\
{\small E-mail: nvchemetov@gmail.com   }%
\\
{\small and }\\
{\small ** CMAF - University of Lisbon, Portugal}}
\date{}
\maketitle
\tableofcontents

\bigskip

\textbf{Abstract.} We consider the flow of an { ideal}
fluid in a 2D-bounded domain, admitting flows through the boundary
of this domain. The flow is described by Euler equations with
\textit{non-homogeneous } Navier slip boundary conditions. These
conditions can be written in the form $\mathbf{v}\cdot
\mathsf{n}=a,$ $2D(\mathbf{v})\mathsf{n}\cdot \mathsf{s}+\alpha \mathbf{v}%
\cdot {\mathsf{s}}=b,$ where the tensor $D(\mathbf{v})$ is the
rate-of-strain of the fluid's velocity $\mathbf{v}$ and\ $\mathsf{(n},%
\mathsf{s)}$ is the pair formed by the normal and tangent vectors
to the boundary. We establish the solvability of this problem in
the class of solutions with $L_{p}-$bounded \ vorticity, $p\in
(2,\infty ].$  To prove the solvability we realize the
passage to the limit in Navier-Stokes equations with vanishing
viscosity.

\bigskip



\newpage

\section{Statement of the problem}

\label{sec0}

The motion of an incompressible inviscous fluid in a domain $\Omega
\subseteq \mathbb{R}^{2}$ is described by the Euler equations
\begin{align}
\mathbf{v}_{t}+\mbox{div}\,(\mathbf{v}\otimes \mathbf{v})-\bigtriangledown
p& =0,\quad \quad (\mathbf{x},t)\in { \Omega _{T}:=\Omega\times (0,T)},  \label{eq1001} \\
\mbox{div}\,\mathbf{v}& =0,\quad \quad (\mathbf{x},t)\in \Omega _{T}
\label{100eq2}
\end{align}%
with a given initial condition
\begin{equation}
\mathbf{v}(\mathbf{x},0)=\mathbf{v}_{0}(\mathbf{x}),\quad \quad \mathbf{x}%
\in \Omega ,  \label{11eq1006}
\end{equation}%
with
\begin{equation}
\mbox{div}\,\mathbf{v}_{0}=0\quad \quad \mbox{ on }\;\Omega .  \label{eq1006}
\end{equation}%
Here $\mathbf{v}(\mathbf{x},t)$ and $p(\mathbf{x},t)$ are the velocity and
the pressure of the fluid at the time $t\in (0,T)$ and the position $\mathbf{%
x}\in \Omega .$

The mathematical theory of the Euler equations was initialized by Gunter
\cite{gunter}, Lichtenstein \cite{lic} and Wolibner
\cite{wol}. They obtained basic results for $%
\Omega =\mathbb{R}^{2}$ and for a bounded domain $\Omega $ of $\mathbb{R}%
^{2} $ with the nonpenetration of $\Gamma :=\partial \Omega $ by the fluid.
The latter nonpenetration condition was written as follows
\begin{equation}
\mathbf{v}\cdot \mathsf{n}=0\quad \quad \mbox{ on }\Gamma _{T}:=\Gamma
\times \lbrack 0,T],  \label{eq3}
\end{equation}%
where $\mathsf{n}$ is an outward normal vector to $\Gamma .$ The global
existence and uniqueness theorems in the 2D-case were obtained by Kato \cite{kato} for the classical solutions  
and  by Yudovich \cite{yu} for the
weak solutions. The problems in domains without
boundaries as $\mathbb{R}^{2},$ $\mathbb{T}^{2}$ (i.e. periodical) and with
the nonpenetration condition \eqref{eq3} have been studied in details.
Discussions about these results can be found by the interested reader in the
books: Antontsev, Kazhikhov, Monakhov \cite{ant}, Lions \cite{lions}, Majda,
Bertozzi \cite{majda} and Temam \cite{temam}.

The aim of this article is to consider Euler equations, admitting a flow of
the fluid through the boundary $\Gamma $ of a bounded domain $\Omega .$ Let
the fluid flow into and out the domain $\Omega $ through the parts $\Gamma
^{-}$ and $\Gamma ^{+}$ of the boundary $\Gamma .$ Assuming that the
quantity of the inflowed and outflowed fluid is equal to $a,$ we have%
\begin{equation}
\mathbf{v}\cdot \mathsf{n}=a\quad \quad \mbox{ on }\Gamma _{T},
\label{eq1.39}
\end{equation}%
such that%
\begin{equation*}
a(\mathbf{x},t):=\left\{
\begin{array}{l}
-|a(\mathbf{x},t)|,\:\:\text{ if }(\mathbf{x},t)\in \Gamma _{T}^{-}:=\Gamma ^{-}\times (0,T); \\
\\
\:\:\:\:\:\:\:\:\:\:\: =0,\:\:\text{ if }(\mathbf{x},t)\in \Gamma _{T}^{0}:=\Gamma ^{0}\times (0,T); \\
\\
\:\:\: |a(\mathbf{x},t)|, \:\:\:\text{ if }(\mathbf{x},t)\in \Gamma _{T}^{+}:=\Gamma ^{+}\times (0,T)%
\end{array}%
\right.
\end{equation*}%
with $\Gamma =\Gamma ^{-}\cup \Gamma ^{0}\cup \Gamma ^{+},$\ $meas(\Gamma
^{+})\neq 0,\,meas(\Gamma ^{-})\neq 0$ and%
\begin{equation}
\int_{\Gamma }a(\mathbf{x},t)\,\,d\mathbf{x}=0\quad \quad \mbox{ for }\;t\in
\lbrack 0,T].  \label{eqC2}
\end{equation}%
The equation \eqref{eq1001} can be written in terms of the vorticity $\omega
:=\mathrm{rot\ }\mathbf{v}\;\;(=\partial _{x_{1}}v_{2}-\partial
_{x_{2}}v_{1})$%
\begin{equation}
\partial _{t}\omega +\mathbf{v}\nabla \omega =0,\quad \quad (\mathbf{x}%
,t)\in \Omega _{T}  \label{eq6}
\end{equation}%
with the initial condition
\begin{equation}
\omega (\mathbf{x},0)=\omega _{0}(\mathbf{x}):=\mathrm{rot\ }\mathbf{v}_{0}(%
\mathbf{x}),\quad \quad \mathbf{x}\in \Omega .  \label{7eq7}
\end{equation}%
Since \eqref{eq6} is a hyperbolic equation, the trajectories of particles
start at the initial moment $t=0$ and on the inflowed boundary $\Gamma ^{-}.$
Therefore, on $\Gamma ^{-}$ we need to put an additional condition. Untill
now three types of this additional condition have been suggested:

\begin{description}
\item[1st type)] The full vector of the velocity $\mathbf{v}$ on $\Gamma
^{-} $ is given
\begin{equation}
\mathbf{v}=\mathbf{b}\quad \quad \text{on }\Gamma _{T}^{-}.  \label{vel}
\end{equation}
\end{description}

The local existence in time of the solution for the problem \eqref{eq1001}-%
\eqref{eq1006}, \eqref{eq1.39}-\eqref{eqC2}, \eqref{vel} was proved by
Kazhikhov \cite{ant}, \cite{kaz}. The solvability of this problem globally
in time is still open.

\begin{description}
\item[2nd type)] The value of the vorticity $\omega $ on $\Gamma ^{-}$ is
given
\begin{equation}
\omega =b\quad \quad \text{on }\Gamma _{T}^{-}.  \label{ome}
\end{equation}
\end{description}

The global existence and uniqueness theorem in the 2D-case of weak solutions
with $L_{p}-$bounded vorticity and $p\in (1,\infty ]$ for the problem %
\eqref{eq1001}-\eqref{eq1006}, \eqref{eq1.39}-\eqref{eqC2}, \eqref{ome} was
obtained by Yudovich \cite{yu}. Later on this problem was considered by
Chemetov, Starovoitov \cite{chem} in the case of \textit{point} sources and
sinks. Here we also mention the article of  Antontsev, Chemetov \cite{ant2}
 where was considered a hyperbolic-elliptic system, 
which describes a flux of superconducting vortices through a domain, for the case of $L_\infty$-bounded
vorticity.

\begin{description}
\item[3d type)] The linear combination of the tangential component of the
viscous stress and the tangential velocity on $\Gamma ^{-}$ is given%
\begin{equation}
2D(\mathbf{v})\mathsf{n}\cdot \mathsf{s}+\alpha \mathbf{v}\cdot {\mathsf{s}}%
=b\quad \quad \mbox{ on }\Gamma _{T}^{-},  \label{eqq1.2}
\end{equation}%
where $\alpha ,b$ are known functions and $\mathsf{s}$ is the tangent to $%
\Gamma $ vector. \ The tensor
\begin{equation*}
D(\mathbf{v}):=\frac{1}{2}[\nabla \mathbf{v}+(\nabla \mathbf{v})^{T}]
\end{equation*}%
is rate-of-strain. We follow the convention that the ordered pair $(\mathsf{n%
},\mathsf{s})$ gives the standard orientation to $\mathbb{R}^{2}$. This
condition is known as \textit{Navier slip boundary condition or friction
condition. }
\end{description}

The inviscid limit of evolutionary Navier--Stokes equations with Navier slip
boundary condition was studied by Clopeau, Mikelic, Robert \cite{clop}. The
existence result for the weak solution with $L_{\infty }-$bounded vorticity
was obtained only for \textit{homogeneous} boundary conditions \eqref{eq3}
and \eqref{eqq1.2} with $b=0.$ In the article of Lopes Filho, Nussenzeng
Lopes, Planas \cite{lopes} this result was generalized in the case of $%
L_{p}- $bounded vorticity with $p\in (2,\infty ),$ but also for homogeneous
boundary conditions \eqref{eq3} and \eqref{eqq1.2} with $b=0.$ Under these
homogeneous boundary conditions, the rate of the convergence of the
Navier-Stokes equations to the Euler equations with respect to the viscosity
parameter was studied by Kelliher \cite{kel}. With the \textit{%
nonhomogeneous }condition\textit{\ }\eqref{eq1.39} and the\textit{\ }Navier
slip boundary condition \eqref{eqq1.2} with $\alpha ,b=0,$ the inviscid
limit of the evolutionary Navier--Stokes equations was investigated by Mucha
\cite{mucha}. Assuming a \textit{geometrical constraint} on the shape of the
domain $\Omega ,$\textit{\ }he proved that the inviscid limit gives the
existence of the solution for the system \eqref{eq1001}-\eqref{eq1006}, %
\eqref{eq1.39}-\eqref{eqC2} with $L_{\infty }-$bounded vorticity. However he
did not prove that the constructed solution satisfies the Navier slip
boundary condition \eqref{eqq1.2} with $\alpha ,b=0.$


This article is a version of the article entitled ”Euler equations
with non-homogeneous Navier slip boundary condition” published in the journal ''Physica D'': Nonlinear Phenomena 237 (2008), 92-105, we refer to \cite{C}. Since the publication of this paper, other articles on this subject have appeared, such as \cite{CC6}-\cite{CC5}, where the theory of boundary layers have been developed for Navier-Stokes and Euler equations with Navier slip boundary condition. A very interesting problem for Euler equations with sources and sinks has been investigated in the article Chemetov, Starovoitov \cite{chem}.  

Moreover a significant progress in the study of the conservative non-linear hyperbolic law coupled with elliptic equations, being a similar system as in the formulation of  Euler equations in terms of vorticity - stream function. This type systems appear in the superconductivity model and Keller-Segel's model. We refer to the articles \cite{ant2, AC2},  \cite{chem1, chem2}, where the Kruzkov approach  has been generalized  using the kinetic method.  We also mention the articles \cite{chem3}-\cite{chem5}, where the kinetic method has been developed to porous media models.  The stochastic case for the non-linear hyperbolic-elliptic systems have been treated in \cite{ACC21}, \cite{CC0}.


In the present paper we consider the problem \eqref{eq1001}-\eqref{eq1006}, %
\eqref{eq1.39}-\eqref{eqC2}, \eqref{eqq1.2} for arbitrary given data $\;a(%
\mathbf{x},t),\;\alpha (\mathbf{x},t),\ b(\mathbf{x},t),\;\mathbf{v}_{0}(%
\mathbf{x}).$ We prove the existence of a weak solution of this problem with
$L_{p}-$bounded vorticity, $p\in (2,\infty ],$ showing that the solution
satisfies the\textit{\ }Navier slip boundary condition \eqref{eqq1.2}.

\textit{The article is divided in the following parts: }

In section \ref{sec444} we give the definition of the weak solution for the
problem \eqref{eq1001}-\eqref{eq1006}, \eqref{eq1.39}-\eqref{eqC2}, %
\eqref{eqq1.2} and formulate the main result of the article. In section \ref%
{sec2} we remember some well-known results, that will be very useful in the
following considerations.

In section \ref{sec44}, we consider a coupled system of the vorticity
equation and a "rot--div" equation on a stream function of the velocity. The
system depends on a \textit{viscosity} and some \textit{regularization
parameters} of the data $a,\alpha ,b,\mathbf{v}_{0}.$ This approximate
problem is solvable, in view of fixed point Shauder%
\'{}%
s argument. In section \ref{sec45}, by a maximum principle, we deduce $%
L_{\infty }$ - estimates for the vorticity and some additional a priori
estimates for the solution of the approximate problem. These estimates allow
us to pass to the limit on the \textit{viscosity} parameter.

In section \ref{secbbb}, we show that\textit{\ solution of the equation %
\eqref{eq6} has to satisfy Gronwall's type inequality} (Theorem \ref%
{lem6sec4234}). As a consequence, in section \ref{sec555}, we derive the
boundedness of the vorticity in $L_{p}$\ with $p\in (2,\infty ],$ uniformly
in time and on the regularization parameters. Then we use a compactness
argument to establish the existence of the weak solution for the original
problem \eqref{eq1001}-\eqref{eq1006}, \eqref{eq1.39}-\eqref{eqC2}, %
\eqref{eqq1.2}.

\vspace{1pt}

\vspace{1pt}

\section{Main and Auxiliary Results}
\setcounter{equation}{0}

\bigskip \bigskip

\subsection{Notations and Main Result}

\label{sec444}

Using the notation, introduced in the books \cite{LadySolonUral68}, \cite%
{LadyUral68}, we shall consider the Sobolev spaces $L_{q}(\Omega ),$ $%
W_{q}^{l}(\Omega ),$ $W_{q}^{l,\,m}(\Omega _{T}),$ $W_{q}^{l}(\Gamma ),$ $%
\,\,\,W_{q}^{l,\,m}(\Gamma _{T})\,$ for $\,\,l,q\geqslant 1$, $m\geqslant 0$
and the H\"{o}lder spaces $C^{l}(\Omega ),$ $C^{l,\,m}(\Omega _{T}),$ $%
C^{l,\,m}(\Gamma _{T})\,\,$ for $\,\,l,m\geqslant 0$, where the index $l$
corresponds to the variable $\mathbf{x}$ and $m$ to the variable $t$ ($%
\,\,l,q,m$ are integer or noninteger).

Let $B$ be a given Banach space. We denote by $C(0,T;B)$ the space of
continuous functions from $[0,T]$ into the Banach space $B$ with the norm
$$
\quad \quad ||u||_{C(0,T;\,\,B)}:=\max_{t\in \lbrack 0,T]}||u(t)||_{B} 
$$  
 and by $L_{q}(0,T;B),$ $1\leqslant q\leqslant \infty $, the space of
measurable functions from $[0,T]$ into the Banach space $B$ with $q$-th
power summable on $[0,T]$, with the norm
\begin{equation*}
||u||_{L_{q}(0,T;\,\,B)}:=(\int_{0}^{T}||u(t)||_{B}^{q}\;dt)^{\frac{1}{q}%
}\quad \quad (:=ess \sup_{0<t<T}||u(t)||_{B}\mbox{ if }q=\infty ).
\end{equation*} 
We introduce also the space 
$V _{2}^{1,0} (\Omega_{T})$  with the norm
$$
||u||_{V _{2}^{1, 0} (\Omega_{T})}=\max_{t \in [0,T] }||u(x,t)||_{L_2 (
\Omega)}+||u_{x}(x,t)||_{L_2 ( \Omega_{T})} .
$$ 
Let us define the following function space%
\begin{equation*}
V_{a}=\{\mathbf{u}\in H^{1}(\Omega ):\;\mbox{div}\mathbf{u}=0\mbox{ in
}\Omega \mbox{ and }\mathbf{u}\cdot \mathsf{n}=a\mbox{ on }\Gamma \}.
\end{equation*}%
From now on let the boundary $\Gamma $ be $C^{2}$-smooth$.$ If we
parametrize the boundary $\Gamma $ using the arc length $s,$ it follows that
$\frac{d\mathsf{n}}{ds}=k{\mathsf{s}},$ where the curvature $k$ of $\Gamma $
is $C$-function, since $\Gamma \in C^{2}.$ The following results are a
direct generalization of Lemma 4.1 and Corollaries 4.2, 4.3 of \cite{kel}.

\begin{lemma}
If $\mathbf{u}\in V_{a}\cap H^{2}(\Omega ),$ then%
\begin{equation*}
(\mathsf{n\cdot }\nabla \mathbf{u)}\cdot \mathsf{s}=\omega (\mathbf{u})-k%
\mathbf{u}\cdot \mathsf{s}+a_{s}^{\prime }\quad \quad \mbox{ and }\quad
\quad (\mathsf{n\cdot }D(\mathbf{u}))\cdot \mathsf{s}=\frac{1}{2}\omega (%
\mathbf{u})-k\mathbf{u}\cdot \mathsf{s}+a_{s}^{\prime },
\end{equation*}%
where $\omega (\mathbf{u})=\mathrm{rot\ }\mathbf{u}$ is the vorticity and $%
a_{s}^{\prime }:=\frac{\partial a}{\partial s}.$ \label{lem4.1}
\end{lemma}

\begin{coro}
\label{coro4.2 copy(1)} A vector $\mathbf{u}\in V_{a}\cap H^{2}(\Omega )$
satisfies the Navier slip boundary condition \eqref{eqq1.2} if and only if
\begin{equation}
\omega (\mathbf{u})=(2k-\alpha )\,\mathbf{u}\cdot \mathsf{s}%
+(b-2\,a_{s}^{\prime }).  \label{(1)}
\end{equation}
\end{coro}

\begin{coro}
\label{coro4.3 copy(1)} For the velocity $\mathbf{u}\in V_{a}\cap C^{2}(%
\overline{\Omega })$, the Navier slip boundary condition \eqref{eqq1.2} with
$\alpha =2k$ is equivalent to
\begin{equation}
\omega (\mathbf{u})=b-2\,a_{s}^{\prime }.  \label{Yudo}
\end{equation}
\end{coro}

That is, any solution of the Euler system \eqref{eq1001}-\eqref{eq1006}, %
\eqref{eq1.39}-\eqref{eqC2}, \eqref{eqq1.2} with $\alpha =2k$ is also a
solution of the Euler system \eqref{eq1001}-\eqref{eq1006}, \eqref{eq1.39}-%
\eqref{eqC2} with Yudovitch%
\'{}%
s type boundary condition \eqref{ome} \ (see \cite{yu}). Let us introduce
the following functions on the boundary $\Gamma _{T}$
\begin{equation}
\left\{
\begin{array}{l}
\gamma :=2k-\alpha ,\;\;g:=b-2\,a_{s}^{\prime },\quad \quad \omega _{\Gamma
}(\mathbf{u}):=\gamma \,\mathbf{u}\cdot {\mathsf{s}}+g; \\
\\
A:=\int_{0}^{s}a(r,t)\,dr\quad \quad \mbox{ for
}t\in \lbrack 0,T],%
\end{array}%
\right.  \label{boundary}
\end{equation}%
where $s$ is the arc length on $\Gamma $. According to (\ref{(1)}), the
system \eqref{eq1001}-\eqref{eq1006}, \eqref{eq1.39}-\eqref{eqC2}, %
\eqref{eqq1.2}\ is equivalent to the problem for the vorticity $\omega $ and
the stream function $h,$ which is a coupling of two systems
\begin{equation}
\left\{
\begin{array}{l}
\omega _{t}+\mbox{div}(\omega \mathbf{v})=0\quad \quad \mbox{ with }\quad
\quad \mathbf{v}=\mathrm{rot}\;h\quad \quad \mbox{ for }(\mathbf{x},t)\in
\Omega _{T}; \\
\\
\omega =\omega _{\Gamma }(\mathbf{v})\quad \quad \mbox{ on }\Gamma
_{T}^{-};\quad \quad \omega (\mathbf{x},0)=\omega _{0}(\mathbf{x})\quad
\quad \mbox{
for }\mathbf{x}\in \Omega%
\end{array}%
\right.  \label{w}
\end{equation}%
and
\begin{equation}
\left\{
\begin{array}{ll}
-\Delta h=\omega & \quad \quad \mbox{ for }(\mathbf{x},t)\in \Omega _{T}; \\
&  \\
h=A(\mathbf{x},t) & \quad \quad \mbox{ for }(\mathbf{x},t)\in \Gamma _{T}.%
\end{array}%
\right.  \label{h}
\end{equation}

\bigskip

\textbf{Regularity of data. }Let $p\in (2,+\infty ].$ \textit{We assume that
the data }$\,$\textit{\thinspace \thinspace }$a,\alpha ,b,\omega _{0}\,\,\,$%
\textit{\ satisfy the conditions} {
\begin{equation}
a \in L_{\infty }(0,T;W_{p}^{1}(\Gamma )), \quad \alpha, b\in
L_{1}(0,T;L_{p}(\Gamma ^{-})),  \quad \omega _{0} \in L_{p}(\Omega
).  \label{eq00sec1}
\end{equation}}
\bigskip

Let us define a weak solution of the Euler equations with the Navier slip
boundary condition.

\begin{definition}
\label{def2sec1} A pair of functions $\{\omega ,h\}$ is said to be a
\underline{weak solution} of the problem \eqref{eq1001}-\eqref{eq1006}, %
\eqref{eq1.39}-\eqref{eqC2}, \eqref{eqq1.2}, if $\omega \in L_{\infty
}(0,T;L_{p}(\Omega )),\;$ $h\in L_{\infty }(0,T;\;W_{q}^{2}(\Omega ))\subset
L_{\infty }(0,T;\;C^{1+\alpha }(\overline{\Omega }))\,$ for $\alpha =1-\frac{%
2}{q}$\ (here $q=p,$\ if $p\in (2,\infty );$ any $q\in (2,\infty )$ if $%
p=\infty )$\ \ and the following equalities hold
\begin{align}
\int_{\Omega _{T}}\omega (\psi _{t}+\mathbf{v}\,\nabla \psi )\;d\mathbf{x}%
dt+\int_{\Omega }\omega _{0}\;\psi (\mathbf{x},0)\;d\mathbf{x}&
=\int_{\Gamma _{T}^{-}}\,a\;\omega _{\Gamma }(\mathbf{v})\,\psi \;d\mathbf{x}%
dt, & &  \label{eq7} \\
-\Delta h=\omega \quad \quad \mbox{ and }\quad \quad \mathbf{v}& =\mathrm{rot%
}\;h\quad \quad \mbox{a.e. in }\Omega _{T}, & &  \label{eq8} \\
\mathbf{v}\cdot \mathsf{n}& =a\quad \quad \mbox{on }\Gamma _{T} & &
\label{eq10}
\end{align}%
for an arbitrary function $\psi \in C^{1,1}(\overline{\Omega }_{T}),$
satisfying the condition%
\begin{equation}
\mbox{ supp}\,(\psi )\subset (\Omega \cup \Gamma ^{-})\times \left[
0,T\right) .  \label{psi1}
\end{equation}%

\end{definition}
{ We shall call the function $\psi $ a \underline{test
function}.}

Our main result is the following theorem.

\begin{teo}
\label{teo2sec1} If the data $\,\,a,\,\,\alpha ,\,\,b,\,\,\omega _{0}$
satisfy (\ref{eq00sec1}), then there exists at least one weak solution $%
\{\omega ,h\}$ of the problem \eqref{eq1001}-\eqref{eq1006}, \eqref{eq1.39}-%
\eqref{eqC2}, \eqref{eqq1.2}.
\end{teo}

\bigskip

\subsection{{ Useful} well-known results}

\label{sec2}

In this section for the convenience of the reader we give well-known results.

According to the potential theory for elliptic equations \cite{vlad}, the
solutions $h_{1},h_{2}$ of the problems
\begin{equation}
\left\{
\begin{array}{ll}
-\Delta h_{1}=F, & \mathbf{x}\in \Omega ; \\
&  \\
h_{1}=0, & \mathbf{x}\in \Gamma%
\end{array}%
\right. \qquad \;\qquad \;\mbox{   and  }\qquad \;\;\qquad \left\{
\begin{array}{ll}
-\Delta h_{2}=0, & \mathbf{x}\in \Omega ; \\
&  \\
h_{2}=f, & \mathbf{x}\in \Gamma%
\end{array}%
\right.  \label{eq13''sec3}
\end{equation}%
can be written in the form
\begin{equation*}
h_{1}(\mathbf{x})=(K_{1}\ast F)(\mathbf{x}):=\int_{\Omega }K_{1}(\mathbf{x},%
\mathbf{y})\,\,F(\mathbf{y})\;d\mathbf{y}\;\mbox{   and  }\;h_{2}=(K_{2}\ast
f)(\mathbf{x}):=\int_{\Gamma }K_{2}(\mathbf{x},\mathbf{y})f(\mathbf{y})\;d%
\mathbf{y.}
\end{equation*}%
Since $\Gamma $ is $C^{2}$-smooth, in the $2-$dimensional case
 { the kernel $%
K_{1}$ satisfies the inequalities%
\begin{equation}
\left\{
\begin{array}{ll}
\left\vert K_{1}(\mathbf{x},\mathbf{y})\right\vert \leqslant C|ln|\mathbf{x}-%
\mathbf{y}||, & \left\vert \nabla _{\mathbf{x}}K_{1}(\mathbf{x},\mathbf{y}%
)\right\vert \leqslant C|\mathbf{x}-\mathbf{y}|^{-1} \\
&  \\
\left\vert \nabla _{\mathbf{x}}\nabla _{\mathbf{y}}K_{1}(\mathbf{x},\mathbf{y%
})\right\vert \leqslant C|\mathbf{x}-\mathbf{y}|^{-2}, & \quad
\mbox{
for any }\mathbf{x},\mathbf{y}\in \Omega .%
\end{array}%
\right.  \label{KK}
\end{equation}%
}

 Taking into
account the potential theory and imbedding theorems of Sobolev
\cite{LadySolonUral68}, \cite{LadyUral68}, \cite{vlad}, we obtain
that the operators $F\rightarrow K_{1}\ast F\;$and $f\rightarrow
K_{2}\ast f$\ \ possess the following properties:

\begin{lemma}
\label{teo2sec2} 1) The function $h_{1}=K_{1}\ast F$ satisfies the following
estimates
\begin{align}
||h_{1}||_{W_{q}^{2}(\Omega )}& \leqslant C||F||_{L_{q}(\Omega )},\quad
\quad \mbox{ if }\quad q\in (1,\infty );  \label{ell-1} \\
\;||h_{1}||_{C^{1+\alpha }(\overline{\Omega })}& \leqslant
C||F||_{L_{q}(\Omega )},\quad \quad \mbox{ if }\quad q\in (2,\infty )\quad %
\mbox{ for }\quad \alpha =1-\frac{2}{q};  \label{ell-5}
\end{align}

2) The function $h_{2}=K_{2}\ast f$ satisfies the following estimates
\begin{align}
||h_{2}||_{W_{q}^{2}(\Omega )}& \leqslant C||f||_{W_{q}^{2-\frac{1}{q}%
}(\Gamma )},\quad \quad \mbox{ if }\quad q\in (1,\infty );  \label{ell-6} \\
||h_{2}||_{C^{1+\alpha }(\overline{\Omega })}& \leqslant C||f||_{W_{q}^{2-%
\frac{1}{q}}(\Gamma )},\quad \quad \mbox{ if }\quad q\in (2,\infty )\quad %
\mbox{ for }\quad \alpha =1-\frac{2}{q}.  \label{ell-8}
\end{align}
\end{lemma}

In the sequel, we shall use the following { extension}
result (see, for instance, \cite{galdi}).

\begin{lemma}
\label{teo33} \ Let $Q\subseteq \mathbb{R}^{n},\quad n\geqslant 1$\ be a
locally lipschitzian domain and let $\;q\in (1,\infty ).$

1) If $u\in W_{q}^{1}(Q),\;$then $u\in W_{q}^{1-\frac{1}{q}}(\partial Q)$ \
and
\begin{equation}
||u||_{W_{q}^{1-\frac{1}{q}}(\partial Q)}\leqslant C||u||_{W_{q}^{1}(Q)};
\label{ell-9}
\end{equation}

2) Any given function $u\in W_{q}^{1-\frac{1}{q}}(\partial Q),$ defined on $%
\partial Q,$ \ can be extended into the domain \ $Q$ in such a way that $%
u\in W_{q}^{1}(Q)$\ \ and
\begin{equation}
||u||_{W_{q}^{1}(Q)}\leqslant C||u||_{W_{q}^{1-\frac{1}{q}}(\partial Q)};
\label{ell-10}
\end{equation}

3) If $u\in W_{q}^{1}(Q),$ then\ for any locally lipschitzian\ domain $%
G\subseteq \mathbb{R}^{n}:Q\subset G,$ there exists an extension of $u$ into
$G,$ such that $u\in W_{q}^{1}(G),$ satisfying $||u||_{W_{q}^{1}(G)}%
\leqslant C||u||_{W_{q}^{1}(Q)}\ \mbox{ and }\ u|_{\partial G}=0.$
\end{lemma}

\section{Construction of approximate solutions}
\setcounter{equation}{0}

\bigskip

The construction of the approximate solutions and the limit transitions on
the viscosity and the regularization parameters will be done for the case
when $p\in (2,\infty ).$ The case $p=\infty $ will be considered just at the
end of the article.\bigskip

\bigskip

\bigskip

\subsection{{ Schauder's} fixed point argument}

\label{sec44} 

\bigskip

Let us define distance functions on $\Gamma .$

\begin{definition}
\label{def2sec copy(1)} Let the distance between any given point $\mathbf{x}%
\in \mathbb{R}^{2}$ and any subset $Q \subseteq \mathbb{R}^{2}$ be defined by
$d(\mathbf{x},Q ):=inf_{\mathbf{y}\in Q }|\mathbf{x}-\mathbf{y}|$. The set of
all points of $\overline{\Omega }$, whose distance to $\Gamma $ ( to $\Gamma
^{-}$ and to $\Gamma ^{+}$) is less than $\sigma $, is denoted by $U_{\sigma
}(\Gamma )$ (respectively, by $U_{\sigma }(\Gamma ^{-})$ and by $U_{\sigma
}(\Gamma ^{+})$). Let $d=d(\mathbf{x})$ be the distance function on $\Gamma $%
, defined by $d(\mathbf{x}):=d(\mathbf{x},\mathbb{R}^{2}\backslash \Omega
)-d(\mathbf{x},\Omega )$ for any$\ \mathbf{x}\in \mathbb{R}^{2}.$
\end{definition}

Since $\Gamma \in C^{2}$, the function $d=d(\mathbf{x})$ has the following
properties

\begin{equation}
d\in C^{2}\quad \mbox{ in  }U_{\sigma _{0}}(\Gamma )\mbox{ for some  }{%
\sigma _{0}}>0\quad \quad \mbox{ and  }\quad \quad \bigtriangledown d=-%
\mathsf{n}\quad \mbox{ on  }\Gamma .  \label{neib}
\end{equation}

 Let us consider a set of smooth functions $\gamma ^{\theta
},$ $g^{\theta },$ $A^{\theta }\in C^{\infty }(\Gamma \times
\left[ 0,T\right] ),$\ $\omega _{0}^{\theta }\in C^{\infty
}(\overline{\Omega }),$ depending on the parameter $\theta \in
(0,\min \left\{ \sigma _{0},\frac{T}{4}\right\} ).$ Let us put
\begin{eqnarray}
a^{\theta
}:=\left( A^{\theta }\right) _{s}^{\prime } \quad \mbox{ on } \Gamma  , \label{AregA}
\end{eqnarray}
here $s$ is
the arc length on $\Gamma$. We can assume that any derivatives
of these
smooth functions are bounded by constants, depending just on a \underline{%
\textit{fixed}} parameter $\theta ,$\ such that
\begin{eqnarray}
\gamma ^{\theta },g^{\theta } &=&0,\quad \mbox{ on the time interval }\left[
0,\theta \right] ;  \notag \\
&&  \notag \\
\omega _{0}^{\theta } &=&0,\quad \mbox{ if }\mathbf{x}\in U_{\theta }(\Gamma
);  \label{reg0}
\end{eqnarray}%
\begin{eqnarray}
&&%
\begin{array}{ll}
\,||\gamma ^{\theta }||_{L_{1}(0,T;L_{p}(\Gamma ^{-}))}\leqslant ||\gamma
||_{L_{1}(0,T;L_{p}(\Gamma ^{-}))}\mathbf{,}\quad \quad & ||g^{\theta
}||_{L_{1}(0,T;L_{p}(\Gamma ^{-}))}\leqslant ||g||_{L_{1}(0,T;L_{p}(\Gamma
^{-}))};%
\end{array}
\notag \\
&&  \notag \\
&&%
\begin{array}{ll}
||\omega _{0}^{\theta }||_{L_{p}(\Omega )}\leqslant ||\omega
_{0}||_{L_{p}(\Omega )}\mathbf{,} & ||A^{\theta }||_{L_{\infty
}(0,T;W_{p}^{2}(\Gamma ))}\leqslant ||A||_{L_{\infty }(0,T;W_{p}^{2}(\Gamma
))},%
\end{array}
\notag \\
&&  \notag \\
&&%
\begin{array}{ll}
||a^{\theta }||_{L_{\infty }(0,T;C(\Gamma ^{-}))}\leqslant ||a||_{L_{\infty
}(0,T;C(\Gamma ^{-}))}\mathbf{,}\quad &
\end{array}
\label{eq122sec2}
\end{eqnarray}%
where the set of functions $\left\{ a^{\theta }\right\} $\ satisfies the
influx and outflux condition on $\Gamma ^{-},\Gamma ^{+}$%
\begin{equation}
a^{\theta }(\mathbf{x},t):=\left\{
\begin{array}{l}
\,-|a^{\theta }(\mathbf{x},t)|,\quad \mbox{ if }\mathbf{x}\in \Gamma ^{-}; \\
\, \quad\quad \quad=0,\quad \mbox{ if }\mathbf{x}\in \Gamma ^{0}; \\
\quad|a^{\theta }(\mathbf{x},t)|,\quad \mbox{ if }\mathbf{x}\in \Gamma ^{+}.%
\end{array}%
\right.  \label{reg1}
\end{equation}%
Moreover, we have the following convergence%
\begin{equation}
\begin{array}{l}
||\gamma ^{\theta }-\gamma ||_{L_{1}(0,T;L_{p}(\Gamma ^{-}))}+||g^{\theta
}-g||_{L_{1}(0,T;L_{p}(\Gamma ^{-}))}\mathop{\longrightarrow}\limits_{{%
\theta }\rightarrow 0}0, \\
\\
||A^{\theta }-A||_{L_{\infty }(0,T;W_{p}^{2}(\Gamma ))}\mathop{%
\longrightarrow}\limits_{{\theta }\rightarrow 0}0, \\
\\
||\omega _{0}^{\theta }-\omega _{0}\,||_{L_{p}(\Omega )}\mathop{%
\longrightarrow}\limits_{{\theta }\rightarrow 0}0.%
\end{array}
\label{eq122sec2-2}
\end{equation}%
From here we { shall} work with these regular initial and boundary conditions $%
\omega _{0}^{\theta },$ $\gamma ^{\theta },$ $g^{\theta },$ $A^{\theta },$ $%
a^{\theta },$ but in the following considerations for the sake of simplicity
of \ notation, we suppress the dependence on ${\theta }$ and write $\omega
_{0},$ $\gamma ,$ $g,$ $A,$ $a,$ respectively. Let us fix positive numbers $%
R>1$ and $\nu \in (0,1).$

Next we construct the pair $\{\omega ,h\}$ as a solution of an auxiliary
\textbf{Problem} $\mathbf{P,}$ which is a coupling of two following systems.

\bigskip

\textbf{Problem} $\mathbf{P.}$ \ \textit{Find $\omega \in W_{2}^{2,1}(\Omega
_{T})$, satisfying the system
\begin{equation}
\left\{
\begin{array}{l}
\omega _{t}+\mbox{div}(\omega \mathbf{v})={\nu }\Delta \omega \quad \quad
\mbox{
and }\quad \quad \mathbf{v}=\mathrm{rot}\;h\quad \quad \mbox{ on }\;\Omega
_{T}; \\
\\
\omega \big|_{\Gamma _{T}}=\omega _{\Gamma }(\mathbf{v});\quad \quad \omega %
\big|_{t=0}=\omega _{0}%
\end{array}%
\right.  \label{eq3sec2}
\end{equation}%
and find $h\in L_{\infty }(0,T;\ W_{2}^{2}(\Omega ))),$ satisfying the
system }%
\begin{equation}
\left\{
\begin{array}{ll}
-\Delta h=\left\langle \omega \right\rangle , & (\mathbf{x},t)\in \Omega
_{T}; \\
&  \\
h=A(\mathbf{x},t), & (\mathbf{x},t)\in \Gamma _{T},%
\end{array}%
\right.  \label{eq4sec2}
\end{equation}%
\textit{where }%
\begin{equation}
\left\langle \omega \right\rangle (\cdot ,t):=\frac{1}{\theta }%
\int_{t}^{t+\theta }\mathit{\ }\left[ \omega (\cdot ,s)\right] _{R}\ ds
\label{R}
\end{equation}%
\textit{with the cut-off function }$\left[ \cdot \right] _{R},$\textit{\
defined as }$\left[ \omega \right] _{R}:=\max \Big\{-R,\;\min \{R,\omega \}%
\Big\}.$\textit{\ We assume that }$\omega =0$\textit{\ outside of the
interval }$\left[ 0,T\right] .$ \bigskip

Of course, the solution $\{\omega ,h\}$ of \textbf{Problem}
$\mathbf{P}$ depends \ on the parameters $R,$ $\nu ,\ \theta ,$
but since in this section the parameters $R,$ $\nu ,\ \theta $ are
fixed, we do not write the dependence of the functions $\omega ,\
h$ on $R,\nu ,\theta $ and shall continue to write $\omega ,$ $h.$
\textit{In the sequel we shall indicate the dependence of
constants and functions on the { parameters} }$R,\nu
,\theta , $\textit{\ if it is necessary}$.$

To prove the solvability of \textbf{Problem} $\mathbf{P},$ we use Schauder's
fixed point argument. Let us introduce the class of functions
\begin{equation}
\Pi :=\{\omega (\mathbf{x},t)\in C(0,T;\;L_{2}(\Omega )):||\omega
||_{C(0,T;\;L_{2}(\Omega ))}\leqslant M\},  \label{eq5sec2}
\end{equation}%
where an exact value of $M$ will be determined below.

\vspace{1pt}First we define the operator $T_{1}$, that transforms a "fixed"
vorticity into a corresponding stream function%
\begin{equation*}
\Pi \ni \widetilde{\omega }\mapsto T_{1}\left( \widetilde{\omega }\right) =h
\end{equation*}%
which is the solution of (\ref{eq4sec2}), where instead of $\omega $ we put
the chosen $\widetilde{\omega }$. Taking into account (\ref{eq13''sec3}),
the function $h$ can be represented in the form%
\begin{equation}
h(\mathbf{x},t)=(K_{1}\ast \left\langle \widetilde{\omega }\right\rangle )(%
\mathbf{x},t)+(K_{2}\ast A)(\mathbf{x},t)  \label{eq66sec2}
\end{equation}%
for a.e. $(\mathbf{x},t)\in \Omega _{T}$ and, by (\ref{ell-5}), (\ref{ell-8}%
) of Lemma \ref{teo2sec2},\ $h$ satisfies the estimate%
\begin{equation}
||h||_{L_{\infty }(0,T;\ C^{1+\alpha }(\overline{\Omega })}\leqslant
C(||\,\left\langle \widetilde{\omega }\right\rangle ||_{L_{ p }(\Omega
_{T})}+||A||_{L_{\infty }(0,T;\;C^{2}(\Gamma ))})\leqslant C,
\label{eq14sec2}
\end{equation}%
for some $\alpha \in (0,1).$ The constants $C=C(R,\theta )$ are \underline{%
independent of $\nu $}, in view of (\ref{eq122sec2}), (\ref{R}).

The second operator $T_{2}$ describes the evolution of the vorticity%
\begin{equation*}
\mathbf{v}=\mathrm{rot}\;h\mapsto T_{2}\left( h\right) =\omega ,
\end{equation*}%
where $\omega $ is the solution of (\ref{eq3sec2}) for the found $h.$ Since $%
\mathbf{v}=\mathrm{rot}\;h\in L_{\infty }(\Omega _{T})$, taking into account
the results of \cite{LadySolonUral68}, Theorem 4.1, p.153 and Theorem 4.2,
p.160, there exists a unique weak solution $\omega $ of (\ref{eq3sec2}),
such that
\begin{align}
||\omega ||_{V_{2}^{1, 0 }(\Omega _{T})}& \leqslant
C_{\ast }, \label{eq9sec2}
\\
||\omega (\cdot ,t_{1})-\omega (\cdot ,t_{2})||_{L_{2}(\Omega )}& \leqslant
\phi (|t_{1}-t_{2}|),\quad \quad \forall t_{1},t_{2}\in \lbrack 0,T]
\label{eq10sec2}
\end{align}
The constant $C_{\ast }$ and the
function $\phi =\phi (t)\in C([0,T]),\;\phi (0)=0$ depend on $\nu
,R,\ \theta $ and the data $\gamma ,\,\,A,\,\,g.$

Setting $M:=C_{\ast }$ in (\ref{eq5sec2}) and $T=T_{2}\circ T_{1}$, by (\ref%
{eq5sec2}), (\ref{eq9sec2}), (\ref{eq10sec2}), we see that $T$ maps the
bounded set $\Pi $ into a compact subset of $\Pi $.

In order to apply Schauder's fixed point theorem, we need to prove that the
operator $T$ is continuous. Let $\widetilde{\omega }_{n},\widetilde{\omega }%
\in \Pi \,\,{}$ and $\widetilde{\omega }_{n}\rightarrow \widetilde{\omega }$
in $C(0,T;\;L_{2}(\Omega ))$. From (\ref{eq66sec2}), (\ref{ell-5}), (\ref%
{ell-8}) of Lemma \ref{teo2sec2} and using that $|\left\langle \widetilde{%
\omega }_{n}\right\rangle |,\ |\left\langle \widetilde{\omega }\right\rangle
|\leqslant R,$ it follows that for any $q>2$
\begin{align}
||\nabla h_{n}-\nabla h||_{C(0,T;\;C^{\alpha }(\overline{\Omega }))}&
\leqslant C||\left\langle \widetilde{\omega }_{n}\right\rangle -\left\langle
\widetilde{\omega }\right\rangle ||_{C(0,T;\;L_{q}(\Omega ))}
\label{eq11sec2} \\
& \leqslant C\cdot (2R)^{\frac{q-2}{q}}||\widetilde{\omega }_{n}-\widetilde{%
\omega }||_{C(0,T;\;L_{2}(\Omega ))}^{2/q}\overset{n\rightarrow \infty }{%
\longrightarrow }0,  \notag
\end{align}%
where $h_{n}$ and $h$ are the solutions of (\ref{eq4sec2}) with $\omega $
replaced by $\widetilde{\omega }_{n}$ and $\widetilde{\omega },$
respectively. Next we consider $\omega _{n}=T\left( \widetilde{\omega }%
_{n}\right) \mbox{
and }\omega =T\left( \widetilde{\omega }\right) .$ Using \eqref{eq11sec2}
and the results of \cite{LadySolonUral68} we conclude that
\begin{equation*}
||\omega _{n}(\cdot ,t)-\omega (\cdot ,t)||_{C(0,T;\;L_{2}(\Omega
))}+||\nabla (\omega _{n}-\omega )||_{L_{2}(0,T;\;L_{2}(\Omega ))}\overset{%
n\rightarrow \infty }{\longrightarrow }0.
\end{equation*}%
Therefore we conclude that the sequence $\omega _{n}$ itself converges to $%
\omega .$ The continuity of the operator $T$ is proved.

We conclude that there exists a fixed point $\omega $ of $T.$ We have shown
the following Lemma.

\vspace{1pt}

\begin{lemma}
\label{lem3sec2} There exists at least one weak solution $\{\omega ,h\}$ of
the systems (\ref{eq3sec2})-(\ref{eq4sec2}), such that for some $\alpha \in
(0,1)$
\begin{equation*}
\omega \in V_{2}^{1,0}(\Omega _{T})\cap C(0,T;\;L_{2}(\Omega )),\quad \quad
h\in C(0,T;\;C^{1+\alpha }(\overline{\Omega })).
\end{equation*}
\end{lemma}

\bigskip

\subsection{Estimates independent of the viscosity}

\bigskip\label{sec45}

In this subsection we derive a priori estimates for the solution $\{\omega ,h\}$
of \textbf{Problem} $\mathbf{P},$ which do not depend on ${\nu \in (0,1).}$
For this reason we shall write $\{{ 
\omega _{\nu},h_{\nu }}\}$ and $\mathbf{%
v_{{\nu }}:=}\mathrm{rot}\;h_{\nu }\mathbf{,}$ having fixed parameters $%
R,\theta .$ 

\textit{In this subsection all constants }$C$\textit{\ }\underline{%
\textit{do not depend on }$\nu $.}

By the theory of parabolic and elliptic equations the constructed solution $%
\{\omega _{{\nu }},h_{{\nu }}\}$ of \textbf{Problem} $\mathbf{P}$ have a
better regularity.

\begin{teo}
\label{teo4sec2} There exists at least one weak solution $\{\omega _{{\nu }},h_{{\nu }}\}$\ of
\textbf{Problem} $\mathbf{P}$, which satisfies the following
\begin{eqnarray*}
h_{{\nu }} &\in &L_{\infty }(0,T;\ W_{p}^{2}(\Omega ))\cap C^{2+\alpha
,\alpha /2}(\overline{\Omega }_{T}),\quad \quad \partial _{t}\,h_{{\nu }}\in
L_{\infty }(0,T;\ W_{p}^{2}(\Omega )) \\
\omega _{{\nu }} &\in &W_{2}^{2,1}(\Omega _{T})\cap C^{\alpha ,\alpha /2}(%
\overline{\Omega }_{T})\quad \quad \mbox{ for some }\alpha \in (0,1),
\end{eqnarray*}%
and%
\begin{align}
||h_{{\nu }}||_{L_{\infty }(0,T;\;C^{1+\alpha }(\overline{\Omega }))}&
\leqslant C,\quad \quad ||h_{\nu }(\cdot ,t)||_{L_{\infty
}(0,T;\;W_{p}^{2}(\Omega ))}\leqslant C,  \label{h1} \\
||\partial _{t}h_{{\nu }}(\cdot ,t)||_{L_{\infty }(0,T;\;W_{p}^{2}(\Omega
))}& \leqslant C,  \label{h2} \\
||\omega _{{\nu }}||_{L_{\infty }(\Omega _{T})}& \leqslant C,  \label{w2}
\end{align}%
with the constants $C=C(R,\theta ),$ which are independent of $\nu .$
\end{teo}

\noindent

\textbf{Proof.} The estimates (\ref{h1}) are a direct consequence of (\ref%
{eq66sec2}), (\ref{eq14sec2}) and (\ref{ell-1}), (\ref{ell-6}) of Lemma \ref%
{teo2sec2}.

Moreover the function $\partial _{t}h_{{\nu }}$ \ \ satisfies the equation
of the system (\ref{eq4sec2}) with $\frac{\partial }{\partial t}\left\langle
\omega _{{\nu }}\right\rangle $ instead of $\left\langle \omega _{{\nu }%
}\right\rangle .$ Since $|\frac{\partial }{\partial t}\left\langle \omega _{{%
\nu }}\right\rangle |\leqslant \frac{2R}{\theta },$ then again applying (\ref%
{ell-1}), (\ref{ell-6}) of Lemma \ref{teo2sec2}, we have (\ref{h2}).

The function $\omega _{{\nu }}$ fulfills the equation%
\begin{equation*}
\partial _{t}\omega _{{\nu }}-{\nu }\Delta \omega _{{\nu }}=F(\mathbf{x}%
,t),\quad \quad (\mathbf{x},t)\in \Omega _{T}
\end{equation*}%
with $F:=-\nabla \omega _{{\nu }}\,\mathbf{v}_{{\nu }}\in L_{2}(\Omega _{T})$
by (\ref{eq9sec2}), (\ref{h1}). Using Theorem 6.1, p.178, \cite%
{LadySolonUral68}, we deduce $\omega _{{\nu }}(\mathbf{x},t)\in
W_{2}^{2,1}(\Omega _{T})$ and also by Theorem 10.1, p.204, \cite%
{LadySolonUral68}, we have $\;\omega _{{\nu }}\in C^{\alpha ,\alpha /2}(%
\overline{\Omega }_{T})$ for some $\alpha \in (0,1).$ Moreover by the theory
of elliptic equations \cite{LadyUral68}, we conclude that $\;\;h_{{\nu }}\in
C^{2+\alpha ,\alpha /2}(\overline{\Omega }_{T}).$

According to the maximum principle for the solution $\omega _{{\nu }}$ of (%
\ref{eq3sec2}) with the help of the boundedness of regular functions $\omega
_{0},\gamma ,g$ \ and (\ref{h1}), we have%
\begin{equation*}
\max_{\overline{\Omega }_{T}}\left\vert \omega _{{\nu }}\right\vert
\leqslant \max \left[ \ \max_{\overline{\Omega }}\ |\omega _{0}(x)|,\quad
||\gamma ||_{L_{\infty }(\Gamma _{T})}\cdot \max_{\Gamma _{T}}|(\mathbf{v}_{{%
\nu }}\cdot {\mathsf{s}})|+||g||_{L_{\infty }(\Gamma _{T})}\right] \leqslant
C,
\end{equation*}%
which proves (\ref{w2}).
 $\ \ \ \ \blacksquare $

\bigskip \vspace{1pt}

\vspace{1pt}

\subsection{Boundary conditions. Limit transition on the viscosity}

\label{sec5}

\vspace{1pt}

In this subsection we consider the viscosity limit $\nu \rightarrow 0$ (the
parameters $R,\theta $ continue to be fixed). Before doing it, we show an
auxiliary lemma, playing an important role in the proof that a limit
function of the sequence $\left\{ \omega _{{\nu }}\right\} $ satisfies the
boundary condition on $\Gamma _{T}^{-}$ in the sense of the equality (\ref%
{eq7}).

\textit{In this subsection all constants }$C$\textit{\ }\underline{%
\textit{do not depend on }$\nu $.}

In view of (\ref{eq122sec2}), (\ref{reg0}), (\ref{h1}), (\ref{h2}) and Lemma %
\ref{teo33}, there exists an { extension}
$\breve{\omega}_{{\nu }}$ of the boundary condition $\omega
_{\Gamma }(\mathbf{v}_{{\nu }})$ and the initial
condition $\omega _{0}$ into the domain $\Omega _{T},$ such that%
\begin{equation}
\left\{
\begin{array}{ll}
\breve{\omega}_{{\nu }}\in W_{p}^{1,1}(\Omega _{T}), & \breve{\omega}_{{\nu }%
}\in L_{\infty }(0,T;\ W_{p}^{2}(\Omega )),\quad \quad \mbox{ such that } \\
||\breve{\omega}_{{\nu }}||_{W_{p}^{1,1}(\Omega _{T})}\leqslant C, & ||%
\breve{\omega}_{{\nu }}||_{L_{\infty }(0,T;\ W_{p}^{2}(\Omega ))}\leqslant
C\quad \quad \mbox{ and } \\
\breve{\omega}_{{\nu }}\big|_{t=0}=\omega _{0}, & \breve{\omega}_{{\nu }}%
\big|_{\Gamma _{T}}=\omega _{\Gamma }(\mathbf{v}_{\nu }).%
\end{array}%
\right.  \label{W3}
\end{equation}
{ Next we use the notations of the Definitions
\ref{def2sec1} and \ref{def2sec copy(1)}.}
\begin{lemma}
\label{lem6sec4} For any \underline{positive} test function $\psi ,$ we have%
\begin{equation}
\lim_{{\sigma }\rightarrow 0}\,\,\left( \,\overline{\lim_{{\nu }\rightarrow
0}}\,\,\frac{1}{{\sigma }}\int_{0}^{T}\int_{[\sigma <d<{2\sigma }]}|\omega _{%
{\nu }}-\breve{\omega}_{{\nu }}|^{p}\;\left( \mathbf{v}_{{\nu }%
}\bigtriangledown d\right) \,\,\psi \;d\mathbf{x}dt\,\right) =0.
\label{eq35sec3}
\end{equation}
\end{lemma}

\noindent

\textbf{Proof.} The function $z_{{\nu }}:=\omega _{{\nu }}-\breve{\omega}_{{%
\nu }}$ satisfies the problem
\begin{equation}
\left\{
\begin{array}{l}
\partial _{t}z_{{\nu }}+\mbox{div}\,(\mathbf{v_{{\nu }}\ }z_{{\nu }})={\nu }%
\Delta z_{{\nu }}+F_{{\nu }},\quad \quad (\mathbf{x},t)\in \Omega _{T}; \\
\\
z_{{\nu }}|_{\Gamma _{T}}=0,\quad \quad \quad \quad z_{{\nu }}|_{t=0}=0,%
\end{array}%
\right.  \label{z}
\end{equation}%
with $F_{{\nu }}:={\nu }\Delta {\breve{\omega}}_{{\nu }}-\partial _{t}{%
\breve{\omega}}_{{\nu }}-\mathbf{v}_{{\nu }}\nabla {\breve{\omega}}_{{\nu }}$%
. By (\ref{h1}), (\ref{W3}), we have%
\begin{equation}
||F_{{\nu }}||_{L_{p}(\Omega _{T})}\leqslant ||\partial _{t}{\breve{\omega}}%
_{{\nu }}||_{L_{p}(\Omega _{T})}+||\nabla {\breve{\omega}}_{{\nu }%
}||_{L_{p}(\Omega _{T})}\,||\mathbf{v}_{{\nu }}||_{C(\overline{\Omega }%
_{T})}+||\Delta {\breve{\omega}}_{{\nu }}||_{L_{\infty }(0,T;\ L_{p}(\Omega
))}\leqslant C.  \label{eq7sec3}
\end{equation}%
Multiplying the equation of (\ref{z}) by $p|z_{{\nu }}|_{\delta
}^{p-1}sgn_{\delta }\left( z_{{\nu }}\right) \,\,\,$ with $|z_{{\nu }%
}|_{\delta }:=\sqrt{z_{{\nu }}^{2}+\delta^{2} } \, $ and $sgn_{\delta }\left( z_{{%
\nu }}\right) :=\frac{z_{{\nu }}}{\sqrt{z_{{\nu }}^{2}+\delta^{2} }},$ we obtain
the inequality%
\begin{align}
\partial _{t}(|z_{{\nu }}|_{\delta }^{p})+\text{div}(\mathbf{v}_{{\nu }%
}\,|z_{{\nu }}|_{\delta }^{p})= \big( {\nu }\,\bigtriangleup |z_{{\nu }}|_{\delta }^{p}+F_{{\nu }}\big) \, p|z_{{%
\nu }}|_{\delta }^{p-1}sgn_{\delta }\left( z_{{\nu }}\right) .
\label{eqAB7sec3}
\end{align}%
Multiplying (\ref{eqAB7sec3}) by an arbitrary \textit{positive} function $%
\eta \in H^{1,1}(\Omega _{T}),$ satisfying the condition (\ref{psi1}), we
derive
\begin{align*}
-\int \int_{\Omega _{T}}\eta _{t}|z_{{\nu }}|_{\delta }^{p}\;d\mathbf{x}dt
-\int \int_{\Omega _{T}}\,\left( \mathbf{v}_{{\nu }}\nabla \eta \right) \
|z_{{\nu }}|_{\delta }^{p}\;d\mathbf{x}dt+\delta ^{p}\int \int_{\Gamma
_{T}}\,\,A\,\eta \;d\mathbf{x}dt-\delta ^{p}\int_{\Omega }\eta (\mathbf{x}%
,0)d\mathbf{x} \\
\leqslant {\nu }\int \int_{\Omega _{T}}|z_{{\nu }}|_{\delta
}^{p}\bigtriangleup \eta \;d\mathbf{x}dt
-\nu \, \delta ^{p}\int \int_{\Gamma
_{T}}\,\,\frac {\partial \eta } {\partial \mathsf{n} } \;d\mathbf{x}dt 
+\int \int_{\Omega _{T}}p|z_{{\nu }%
}|_{\delta }^{p-1}|F_{{\nu }}|\,\,\eta \;d\mathbf{x}dt.
\end{align*}%
Taking $\delta \rightarrow 0,$ we obtain%
\begin{align*}
& -\int \int_{\Omega _{T}}\,\left( \mathbf{v}_{{\nu }}\nabla \eta \right) \
|z_{{\nu }}|^{p}\;d\mathbf{x}dt \\
& \leqslant \int \int_{\Omega _{T}}\eta _{t}|z_{{\nu }}|^{p}\;d\mathbf{x}dt+{%
\nu }\int \int_{\Omega _{T}}|z_{{\nu }}|^{p}|\bigtriangleup \eta |\;d\mathbf{%
x}dt+\int \int_{\Omega _{T}}p|z_{{\nu }}|^{p-1}|F_{{\nu }}|\,\,\eta \;d%
\mathbf{x}dt.
\end{align*}%
The functions $z_{{\nu }}$ are uniformly bounded in $L_{\infty }(\Omega
_{T}) $ by a constant $C=C(R,\theta ),$ which is independent of $\nu $ by
the estimates (\ref{w2}) and (\ref{W3}). \ Therefore%
\begin{equation}
-\int \int_{\Omega _{T}}\left( \mathbf{v}_{{\nu }}\nabla \eta \right)
\;|z_{\nu }|^{p}\;d\mathbf{x}dt\leqslant C\int \int_{\Omega _{T}}(|\eta
_{t}|+\nu |\bigtriangleup \eta |+|F_{{\nu }}|\,\eta )\ d\mathbf{x}dt.
\label{eq36sec3}
\end{equation}%
Let us take an arbitrary test function $\psi $ and put the function $\eta :=(1-{%
\mathbf{1}_{\sigma }(\mathbf{x})})\,\psi $ in (\ref{eq36sec3}) with%
\begin{equation*}
{\mathbf{1}_{\sigma }}(\mathbf{x}):=\left\{
\begin{array}{l}
0,\quad \mbox{ if }d(\mathbf{x)\in }[0,\sigma ); \\
\\
\frac{d-\sigma }{\sigma },\quad \mbox{ if }d(\mathbf{x)\in }[\sigma ,2\sigma
); \\
\\
1,\quad \mbox{ if }\vec{x}\in \Omega \backslash U_{2\sigma }(\Gamma ).%
\end{array}%
\right.
\end{equation*}%
Then we deduce%
\begin{align}
\int \int_{\Omega _{T}\cap \lbrack \sigma <d<{2\sigma }]}\frac{\left(
\mathbf{v}_{{\nu }}\nabla d\right) }{{\sigma }}\;\psi \,|z_{{\nu }}|^{p}\,d%
\mathbf{x}dt& \leqslant C\int \int_{\Omega _{T}}\left\{ \nu |\bigtriangleup
\eta |+\right.  \notag \\
& \left. +(1-{\mathbf{1}_{\sigma }})\left( |\psi _{t}|+|F_{\nu }|\psi
+|\nabla \psi |\right) \right\} \;d\mathbf{x}dt.  \label{eq1122sec2}
\end{align}%
By (\ref{h1}) the set $\{\mathbf{v}_{{\nu }}(\cdot ,t)\}$ is uniformly
continuous on $\overline{\Omega },$ independently of ${\nu }\in \left(
0,1\right) $ and $t\in \lbrack 0,T].$ With the help of (\ref{neib}), (\ref%
{reg1}), we get that there exists some $\sigma _{1}<\sigma _{0},$
independent of $\nu \in (0,1)$ and $t\in \lbrack 0,T],$ such that
\begin{equation*}
\left\{ \,\mathbf{v}_{{\nu }}\bigtriangledown d\right\} \,(\mathbf{x}%
,t)=\left\{
\begin{array}{ll}
>0, & \mbox{ if }(\mathbf{x},t)\in U_{\sigma _{1}}(\Gamma ^{-})\times
\lbrack 0,T]; \\
&  \\
<0, & \mbox{ if }(\mathbf{x},t)\in U_{\sigma _{1}}(\Gamma ^{+})\times
\lbrack 0,T].%
\end{array}%
\right.
\end{equation*}%
For the test function $\psi ,$ {  the} condition
(\ref{psi1}) means that there exists $\overline{\sigma },$ such
that
\begin{equation}
\psi (\mathbf{x},t)=0\quad \quad \mbox{ for  }\left\{
\begin{array}{l}
\,\,\forall \mathbf{x}\in \overline{\Omega }\,\mbox{ and }\,\,\forall t\in
\lbrack T-\overline{\sigma },T]; \\
\,\forall \mathbf{x}\in U_{\overline{\sigma }}\,(\Gamma ^{+}\cup \Gamma
^{0})\,\mbox{ and }\,\,\forall t\in \lbrack 0,T].%
\end{array}%
\right.  \label{ff}
\end{equation}%
Therefore for any $\sigma \in (0,\ \frac{\min \left( \sigma _{1},\overline{%
\sigma }\right) }{2}),$ according (\ref{eq1122sec2}), we have%
\begin{align*}
0& \leqslant \overline{\mathop{\lim}\limits_{{\nu }\rightarrow 0}}\,\,\frac{1%
}{\sigma }\int_{0}^{T}\int_{[\sigma <d<{2\sigma }]}\left( \mathbf{v}_{{\nu }%
}\nabla d\right) \;\psi \,|z_{{\nu }}|^{p}\,d\mathbf{x}dt \\
& \leqslant C\int \int_{\Omega _{T}}(1-{\mathbf{1}_{\sigma }})\left( |\psi
_{t}|+\psi |F_{{\nu }}|+|\nabla \psi |\right) \;d\mathbf{x}dt\;d\mathbf{x}.
\end{align*}%
This implies the property (\ref{eq35sec3}), in view of (\ref{eq7sec3}) and
since $(1-{\mathbf{1}_{\sigma }(\mathbf{x})})\underset{{\sigma }\rightarrow 0%
}{\longrightarrow }0$ in $\Omega _{T}.$ $\ \ \blacksquare $

\bigskip

From (\ref{h1})-(\ref{w2}) we conclude that there exists a subsequence of $%
\{\omega _{{\nu }},h_{{\nu }}\},$ such that
\begin{align}
h_{{\nu }}& \rightharpoonup h,\quad \partial _{t}h_{{\nu }}\rightharpoonup \
\partial _{t}h\quad \quad \mbox{ weakly}-\ast \mbox{ in }L_{\infty
}(0,T;\,W_{p}^{2}(\Omega )),  \notag \\
\omega _{{\nu }}& \rightharpoonup \omega \quad \quad \mbox{ weakly}-\ast
\mbox{ in
}L_{\infty }(\Omega _{T}),  \label{eq31sec3}
\end{align}%
that also implies%
\begin{equation}
\mathbf{v}_{{\nu }}\rightarrow \mathbf{v:=}\mathrm{rot}\;h\quad \quad
\mbox{   in
}L_{p}(\Omega _{T}).  \label{4.9}
\end{equation}%
The limit functions $\{\omega ,h\}$ fulfills the estimates (\ref{h1})-(\ref%
{w2}).

By (\ref{eq31sec3}) and the representation (\ref{eq66sec2}), it is easy to
check that the function $h$ satisfies the system (\ref{eq4sec2}), such that
\begin{equation}
h\in L_{\infty }(0,T;\,C^{1+\alpha }(\overline{\Omega }))\quad \quad
\mbox{
for  }\alpha =1-\frac{2}{p},  \label{Bpresent}
\end{equation}%
\begin{equation}
h=h_{1}+h_{2}\quad \quad \mbox{
with   }\quad \quad h_{1}:=K_{1}\ast \left\langle \omega \right\rangle
,\,\,\,h_{2}:=K_{2}\ast A.  \label{Apresent}
\end{equation}%
Let us prove that the pair $\{\omega ,\mathbf{v}\}$ satisfies (\ref{eq7}).
Multiplying the equation of (\ref{eq3sec2}) by $\eta _{\sigma }:={\mathbf{1}%
_{\sigma }}\,\psi ,$ where $\psi $ is a test function and integrating it
over $\Omega _{T},$ we derive
\begin{align}
0=\{\int_{\Omega _{T}}[\omega _{{\nu }}(\psi _{t}& +\mathbf{v}_{{\nu }%
}\nabla \psi )]{\mathbf{1}_{\sigma }}+{\nu \ }\omega _{{\nu }}\bigtriangleup
\eta _{\sigma }\;d\mathbf{x}dt+\int_{\Omega }\omega _{0}(\mathbf{x})\,\eta
_{\sigma }(\mathbf{x},0)\;d\mathbf{x}\,\}+ & &  \notag  \label{ddsec3} \\
& +\frac{1}{{\sigma }}\int_{0}^{T}\int_{[\sigma <d<{2\sigma }]}\omega _{{\nu
}}\,(\mathbf{v}_{{\nu }}\nabla d)\,\psi \;d\mathbf{x}dt=I^{{\nu },\sigma
}+J^{{\nu },\sigma }. & &
\end{align}%
Using (\ref{eq31sec3}), (\ref{4.9}) and ${\mathbf{1}_{\sigma }}%
\mathop{\longrightarrow}\limits_{{\sigma }\rightarrow 0}1$ in $\Omega _{T}$
and $\Omega ,$ we have%
\begin{equation*}
\lim_{{\sigma }\rightarrow 0}\,\,\left( \,\lim_{{\nu }\rightarrow 0}I^{{\nu }%
,\sigma }\right) =\int_{\Omega _{T}}\omega (\psi _{t}+\mathbf{v}\nabla \psi
)\;d\mathbf{x}dt+\int_{\Omega }\omega _{0}\psi (\mathbf{x},0)\;d\mathbf{x.}
\end{equation*}

Moreover%
\begin{align*}
J^{\nu ,\sigma }=[\frac{1}{\sigma }\int_{0}^{T}\int_{[\sigma <d<2\sigma ]}z_{%
{\nu }}\,(\vec{v}_{{\nu }}\nabla d)\,\psi \;d\mathbf{x}dt\}]+& \\
+[\frac{1}{\sigma }\int_{0}^{T}\int_{[\sigma <d<2\sigma ]}& \breve{\omega}_{{%
\nu }}\,(\mathbf{v}_{{\nu }}\nabla d)\,\psi \;d\mathbf{x}dt\}]=J_{1}^{\nu
,\sigma }+J_{2}^{\nu ,\sigma }.
\end{align*}%
According to Lemma {\ref{lem6sec4}}, we obtain
\begin{equation*}
\lim_{{\sigma }\rightarrow 0}\,\,\left( \,\overline{\lim_{{\nu }\rightarrow
0}}|J_{1}^{{\nu },\sigma }|\right) =0.
\end{equation*}%
By (\ref{h1}), (\ref{h2}) the set of functions $\mathbf{v}_{{\nu }%
}\bigtriangledown d$ is uniformly continuous on $\overline{\Omega }$ for all
$t\in \lbrack 0,T]$, independently of ${\nu }$ and the trace of $\mathbf{v}_{%
{\nu }}\bigtriangledown d$ on $\Gamma _{T}$ is equal to $-a$ for every point
$(\mathbf{x},t)\in \Gamma _{T}.$ By (\ref{W3}) the function $\breve{\omega}_{%
{\nu }}$ has the trace $\omega _{{\Gamma }}(\mathbf{v}_{\nu })$\ on the
boundary $\Gamma _{T}^{-}$ \ and, in view of (\ref{eq31sec3}), we have the
convergence
\begin{equation*}
\omega _{{\Gamma }}(\mathbf{v}_{\nu })\underset{{\nu }\rightarrow 0}{%
\longrightarrow }\omega _{{\Gamma }}(\mathbf{v})\quad \quad \mbox{ in }\quad
L_{1}(\Gamma _{T}^{-}),
\end{equation*}%
that gives%
\begin{equation*}
\lim_{{\sigma }\rightarrow 0}\,\,(\lim_{{\nu }\rightarrow 0}\,\,J_{2}^{{\
\nu },\sigma })=-\int_{0}^{T}\int_{\Gamma _{T}^{-}}\omega _{{\Gamma }}(%
\mathbf{v})\,a\,\psi \;d\mathbf{x}dt.
\end{equation*}%
Therefore the pair $\{\omega ,\,h\}$ satisfies the equation (\ref{eq7}).

\vspace{1pt}

\section{Existence result}
\setcounter{equation}{0}

\subsection{Gronwall's type inequality for the vorticity}

\label{secbbb}

\bigskip

In the previous section we have shown the existence of the
functions $\omega ,h,$ satisfying the equalities (\ref{eq7}), (\ref{eq4sec2}%
) with $\mathbf{v=}\mathrm{rot}\;h$ and the estimates (\ref{h1})-(\ref{w2})$%
. $ In this subsection we continue to have fixed parameters $R, \theta$, denoting by $C$ constants
which may depend on $R, \theta$.

The main objective is to show \ that the pair $\left\{
\omega ,\mathbf{v}\right\} $ fulfills the following Gronwall's type inequality. 

\begin{teo}
\label{lem6sec4234} For all$\;t_{0}\in \lbrack 0,T],$ we have
\begin{equation}
\int_{\Omega }|\omega (\mathbf{x},t_{0})|^{p}\;d\mathbf{x}\leqslant
\int_{\Omega }|\omega _{0}|^{p}\;d\mathbf{x}+\int_{0}^{t_{0}}\int_{\Gamma
^{-}}a\;|\omega _{\Gamma }(\mathbf{v})|^{p}\;d\mathbf{x}dt.  \label{A}
\end{equation}
\end{teo}

\textbf{\vspace{1pt}}

\vspace{1pt}

\begin{remark}
In\ fact we may \ correctly define a concept of the trace of $\omega $\ on
the boundary $\Gamma ^{+}$\ and show a more strong equality%
\begin{equation*}
\int_{\Omega }|\omega (\mathbf{x},t_{0})|^{p}\;d\mathbf{x+}%
\int_{0}^{t_{0}}\int_{\Gamma ^{+}}\;a\;|\omega |^{p}\;d\mathbf{x}%
dt=\int_{\Omega }|\omega _{0}|^{p}\;d\mathbf{x}+\int_{0}^{t_{0}}\int_{\Gamma
^{-}}\;a\;|\omega _{\Gamma }(\mathbf{v})|^{p}\;d\mathbf{x}dt.
\end{equation*}
\end{remark}

Before proving this theorem, let us show a lemma, which is an adaptation of
Lemma \ref{lem6sec4} to the equation \ (\ref{eq7}), which is of hyperbolic
type.

Put $\Omega ^{\theta }:=\left\{ \mathbf{x}\in \mathbb{R}^{2}:\ d(\mathbf{x}%
,\Omega \mathbf{)}<\theta \right\} $ and $\Omega _{T}^{\theta }:=\Omega
^{\theta }\times \left[ -\theta ,T+\theta \right] $. In view of Lemma \ref%
{teo33} and (\ref{reg0}), (\ref{eq122sec2}), since $\mathbf{v}\in
W_{p}^{1,1}(\Omega _{T}),$ there exists an { extension}
$\breve{\omega}$ of the boundary conditions $\omega _{\Gamma
}(\mathbf{v}):=\gamma \,\mathbf{v}\cdot
{\mathsf{s}}+g$ and the initial conditions $\omega _{0}$ into the domain $%
\Omega _{T}^{\theta },$ such that
\begin{equation}
\left\{
\begin{array}{ll}
\breve{\omega}\in W_{p}^{1,1}(\Omega _{T}^{\theta }),\quad \quad
\mbox{ with
} & ||\breve{\omega}||_{W_{p}^{1,1}(\Omega _{T})}\leqslant C\quad \quad
\mbox{
and } \\
&  \\
\breve{\omega}\big|_{t=0}=\omega _{0}, & \breve{\omega}\big|_{\Gamma
_{T}}=\omega _{\Gamma }(\mathbf{v}),%
\end{array}%
\right.  \label{W}
\end{equation}%
\bigskip where the constant $C$ depends on $\theta ,R$ and $||\mathbf{v||}%
_{W_{p}^{1,1}(\Omega _{T})}.$

Let $b\in L_{q}\,(0,T;\;L_{r}(\Omega ))$ for some $q,r\in \lbrack
\,1,+\infty )$ be a given function on $\Omega _{T}$. We can assume that the
function $b$ is continued outside of $\Omega _{T}$ by zero value. In the
sequel we regularize the function $b$ by the convolution $b^{\varepsilon
}:=\rho _{\varepsilon }\ast b$ with $\rho _{\varepsilon }(\mathbf{x},t):=%
\frac{1}{\varepsilon ^{3}}\rho (\frac{x_{1}}{\varepsilon })\,\rho (\frac{%
x_{2}}{\varepsilon })\,\rho (\frac{t}{\varepsilon }).$ The kernel $\rho (\xi
)$ is equal to $0,$ if $|\xi |\geq 1;$ $K\exp \Bigl(-\frac{1}{1-|\xi |^{2}}%
\Bigr),$ if $|\xi |<1,$ where the constant $K$ is chosen such that $\int_{%
\mathbb{R}}\rho (\xi )\,d\xi =1.$ From the properties of the regularization
operator $(\cdot )^{\varepsilon }$
\begin{equation}
b^{\varepsilon }\underset{{\varepsilon }\rightarrow 0}{\longrightarrow }%
b\quad \quad \mbox{ in
}\quad \quad L_{q}\,(\left[ -\theta ,T+\theta \right] ,L_{r}\,(\Omega
^{\theta })).  \label{(B)}
\end{equation}

\begin{lemma}
\label{lem6sec44} For any \underline{positive} test function $\psi ,$ we
have
\begin{align}
\lim_{{\sigma }\rightarrow 0}\,\,\left( \,\,\frac{1}{{\sigma }}%
\int_{0}^{T}\int_{[\sigma <d<{2\sigma }]}|\omega -\breve{\omega}|^{p}\,\,%
\mathbf{v}\bigtriangledown d\,\,\psi \;d\mathbf{x}dt\,\right) & =0,
\label{bo} \\
\lim_{{\sigma }\rightarrow 0}\,\,\left( \frac{1}{{\sigma }}\int_{0}^{{\sigma
}}\int_{\Omega }|\omega -\breve{\omega}|^{p}\,\,\psi \,d\mathbf{x}dt\right)
& =0.  \label{to}
\end{align}
\end{lemma}

\noindent

\textbf{Proof. } The function $z:=\omega -\breve{\omega}$ satisfies the
following equality
\begin{equation*}
\int_{\Omega _{T}}z(\eta _{t}+\mathbf{v}\,\nabla \eta )\;d\mathbf{x}%
dt=\int_{\Omega _{T}}F\eta \;d\mathbf{x}dt\quad \quad \mbox{ with }\quad
\quad F:=\partial _{t}{\breve{\omega}}+\mbox
{div}\,(\mathbf{v}\,{\breve{\omega}}),
\end{equation*}%
for any test function $\eta $. By (\ref{Bpresent}), (\ref{W}), we have
\begin{equation}
||F||_{L_{p}(\Omega _{T})}\leqslant ||\partial _{t}{\breve{\omega}}%
||_{L_{p}(\Omega _{T})}+||\nabla {\breve{\omega}}||_{L_{p}(\Omega _{T})}\,||%
\mathbf{v}||_{L_{\infty }(0,T;C(\overline{\Omega }))}\leqslant C,
\label{est}
\end{equation}%
where the constant $C$ depends on $\theta ,R,$ $||\mathbf{v||}%
_{W_{p}^{1,1}(\Omega _{T})}$ and $||\mathbf{v||}_{L_{\infty
}(0,T;W_{p}^{1}(\Omega ))}.$ Let us assume that the functions $z,\;F$ are
continued outside of the domain $\overline{\Omega }_{T}$ by zero and the
velocity $\mathbf{v}$ is continued into the domain $\Omega ^{\theta
} \times \left[ -\theta ,T+\theta \right] $ by some function%
\begin{equation}
\widehat{\mathbf{v}}\in L_{\infty }(\left[ -\theta ,T+\theta \right]
;\;W_{p}^{1}(\Omega ^{\theta })),\quad \quad \mbox{ such that }\quad \quad %
\mbox{div}\,\widehat{\mathbf{v}}=0 .  \label{cont}
\end{equation}%
From the properties of the regularization operator $(\cdot )^{\varepsilon }$
(see Theorem II.1 \cite{dip}), we see that the function $z^{\varepsilon
}:=(\omega -\breve{\omega})^{\varepsilon }$ satisfies the equality
\begin{equation*}
\int_{\Omega _{T}^{\theta }}\left( z_{t}^{\varepsilon }+\widehat{\mathbf{v}}%
\nabla z^{\varepsilon }\right) \eta \;d\mathbf{x}dt=\int_{\Omega
_{T}^{\theta }}\left( r^{\varepsilon }-F^{\varepsilon }\right) \eta \;d%
\mathbf{x}dt,
\end{equation*}%
with
\begin{equation}
r^{\varepsilon }:=\mbox
{div}\,\left[ \widehat{\mathbf{v}}\,{z}^{\varepsilon }-\,(\widehat{\mathbf{v}%
}\,{z})^{\varepsilon }\right] \;\longrightarrow 0\quad \quad \mbox
{ in }\quad L_{r}(\left[ -\theta ,T+\theta \right] ;\;W_{p}^{1}(\Omega
^{\theta })),  \label{BB}
\end{equation}%
$\forall \,\,r\in \left[ 1,+\infty \right) $ and for any function $\eta (%
\mathbf{x},t)\in C_{0}^{1,1}(\Omega _{T}^{\theta }),$ satisfying the
condition (\ref{psi1}). From the previous, we easily deduce the equality
\begin{equation*}
-\int \int_{\Omega _{T}^{\theta }}|z^{\varepsilon }|^{p}(\eta _{t}+\widehat{%
\mathbf{v}}\,\nabla \eta )\;d\mathbf{x}dt=\int \int_{\Omega _{T}^{\theta
}}(r^{\varepsilon }-F^{\varepsilon })\ p|z^{\varepsilon
}|^{p-1}sgn(z^{\varepsilon })\ \eta \;d\mathbf{x}dt.
\end{equation*}%
Since $z=0,$ $\ F=0$ outside of the domain $\Omega _{T}$ and $\widehat{%
\mathbf{v}}=\mathbf{v}$\ in $\Omega _{T},$ with the help of (\ref{w2}), \ (%
\ref{W}), (\ref{(B)}), \ (\ref{est})-(\ref{BB}), the limit transition on $%
\varepsilon \rightarrow 0$ gives that $|z|^{p}$ fulfills the equality
\begin{equation}
-\int \int_{\Omega _{T}}|z|^{p}\ (\eta _{t}+\mathbf{v}\,\nabla \eta )\;d%
\mathbf{x}dt=\int \int_{\Omega _{T}}-F\ p|z|^{p-1}sgn(z)\ \eta \;d\mathbf{x}%
dt.  \label{lut}
\end{equation}%
for any test function $\eta .$

In view of the regularity (\ref{Bpresent}), the function $\mathbf{v}(\cdot
,t)$ is uniformly continuous on $\overline{\Omega }$, independently of $t\in
\lbrack 0,T]$. Hence, taking into account (\ref{reg1}) and (\ref{neib}),
there exists $\sigma _{2}<\sigma _{0},$ such that
\begin{equation}
\left\{ \mathbf{v}\bigtriangledown d\right\} \,(\mathbf{x},t)=\left\{
\begin{array}{ll}
>0, & \mbox{ if }(\mathbf{x},t)\in U_{\sigma _{2}}(\Gamma ^{-})\times
\lbrack 0,T]; \\
&  \\
<0, & \mbox{ if }(\mathbf{x},t)\in U_{\sigma _{2}}(\Gamma ^{+}\cup \Gamma
^{0})\times \lbrack 0,T].%
\end{array}%
\right.  \label{fff}
\end{equation}%
Let us fix a \textit{positive} test function $\psi $.\ Accounting (\ref{ff})
and choosing in the equality (\ref{lut})

\vspace{1pt}a) $\eta (\mathbf{x,t}):=(1-{\mathbf{1}_{\sigma }(}\mathbf{x)}%
)\cdot \psi (\mathbf{x,t})$, \ \ we derive%
\begin{eqnarray*}
0 &\leqslant &\frac{1}{\sigma }\int_{0}^{T}\int_{[\sigma <d<{2\sigma }%
]}|z|^{p}\,\mathbf{v}\nabla d\;\psi \;d\mathbf{x}dt \\
&\leqslant &C\int \int_{\Omega _{T}}(1-{\mathbf{1}_{\sigma }(}\mathbf{x)}%
)\cdot \left\{ |F|\;\psi |z|^{p-1}+\left[ |\partial _{t}\psi |+|\nabla \psi |%
\right] \;|z|^{p}\right\} \;d\mathbf{x}dt
\end{eqnarray*}%
for any fixed $\sigma \in (0,\frac{\min \left( \sigma _{2},\overline{\sigma }%
\right) }{2})$ (this inequality implies the property \eqref{bo}, in view of (%
\ref{fff}) and the convergence $(1-{\mathbf{1}_{\sigma }(}\mathbf{x)})\,%
\underset{\sigma \rightarrow 0}{\longrightarrow }0$ in $\Omega _{T});$

b) \ \ $\eta (\mathbf{x,t}):=(1-1{_{\sigma }(t}\mathbf{)})\cdot \psi (%
\mathbf{x,t})$ with%
\begin{equation*}
1{_{\sigma }(t)}:=\left\{
\begin{array}{l}
0,\quad \mbox{ if }t\in (-\infty ,\sigma ); \\
\\
\frac{t-\sigma }{\sigma },\quad \mbox{ if }t\in \lbrack \sigma ,2\sigma );
\\
\\
1,\quad \mbox{ if }t\in \lbrack 2\sigma ,+\infty ),%
\end{array}%
\right.
\end{equation*}%
we derive
\begin{eqnarray*}
0 &\leqslant &\,\frac{1}{\sigma }\int_{{\sigma }}^{{2\sigma }}\int_{\Omega
}\,|z|^{p}\,\;\psi \;d\mathbf{x}dt \\
&\leqslant &C\int \int_{\Omega _{T}}(1-1{_{\sigma }(t}\mathbf{)})\cdot
\left\{ |F|\;\psi |z|^{p-1}+\left[ |\partial _{t}\psi |+|\nabla \psi |\right]
\;|z|^{p}\right\} \;d\mathbf{x}dt
\end{eqnarray*}%
for any fixed $\sigma \in (0,\frac{\min \left( \sigma _{2},\overline{\sigma }%
\right) }{2}).$ It implies the property \eqref{to}, since $(1-1{_{\sigma }(t}%
\mathbf{)})\,\underset{\sigma \rightarrow 0}{\longrightarrow }0$ in $\Omega
_{T}$. This conclude the proof of Lemma. $\;\;\;$ $\;\;\;\;\blacksquare $

\bigskip

Now we are able to prove Theorem \ref{lem6sec4234}. \ Since $\left\{ \omega
,h,\mathbf{v}\right\} $ satisfies the equation (\ref{eq7}) and the estimates
(\ref{h1})-(\ref{w2}), using the theory of \cite{dip}, we have that the pair
$\left\{ |\omega |^{p},\mathbf{v}\right\} $ fulfills the equality
\begin{equation}
\int_{\Omega _{T}}|\omega |^{p}\left( \eta _{t}+\mathbf{v\cdot }\nabla \eta
\right) \;d\mathbf{x}dt=0\quad \mbox{ for any
}\quad \eta \in C_{0}^{1,1}(\Omega _{T}).  \label{horosho}
\end{equation}%
If we put in this equality $\eta (\mathbf{x},t):=\mathbf{1}_{\sigma }(%
\mathbf{x})\cdot \phi (t),\ $ where $\sigma \in (0,\frac{\sigma _{0}}{2})$
and $\phi \in C_{0}^{1}([0,T]).\ $ Then we see that the integral
\begin{equation}
\int_{\Omega }\;|\omega (\mathbf{x},t)|^{p}\ \mathbf{1}_{\sigma }(\mathbf{x}%
)\;d\mathbf{x}\in W_{1}^{1}([0,T])\subset C([0,T]).  \label{hor}
\end{equation}%
Moreover, choosing $\eta (\mathbf{x},t):=\left[ 1_{\varepsilon
}(t)-1_{\varepsilon }(t+\varepsilon -t_{0})\right] \cdot \mathbf{1}_{\sigma
}(\mathbf{x})\ $ for any $\varepsilon ,\sigma \in (0,\frac{\min \left\{
\sigma _{0},T\right\} }{4})$ and { arbitrary} fixed $t_{0}\in \lbrack 0,T]$ in (%
\ref{horosho}), we derive
\begin{eqnarray*}
-\left[ \frac{1}{\varepsilon }\int_{\varepsilon }^{2\varepsilon
}\int_{\Omega }\;|\omega |^{p}\mathbf{1}_{\sigma }(\mathbf{x})\;d\mathbf{x}%
dt-\frac{1}{\varepsilon }\int_{t_{0}}^{t_{0}+\varepsilon }\int_{\Omega
}\;|\omega |^{p}\mathbf{1}_{\sigma }(\mathbf{x})\;d\mathbf{x}dt\right] && \\
-\frac{1}{\sigma }\int_{0}^{T}\int_{[\sigma <d<2\sigma ]}\;\left[
1_{\varepsilon }(t)-1_{\varepsilon }(t+\varepsilon -t_{0})\right] | &&\omega
|^{p}\;\mathbf{v}\nabla d\;d\mathbf{x}dt=0.
\end{eqnarray*}%
Taking the limit on $\varepsilon \rightarrow 0,$ by (\ref{hor}), (\ref{to})
and (\ref{W}), we obtain
\begin{equation*}
\int_{\Omega }\;|\omega (\mathbf{x},t_{0})|^{p}\mathbf{1}_{\sigma }(\mathbf{x%
})\;d\mathbf{x}\leqslant \int_{\Omega }\;|\omega _{0}(\mathbf{x})|^{p}%
\mathbf{1}_{\sigma }(\mathbf{x})\;d\mathbf{x}dt+\frac{1}{\sigma }%
\int_{0}^{t_{0}}\int_{U_{\sigma }(\Gamma ^{-})}|\omega |^{p}\;\mathbf{v}%
\nabla d\;d\mathbf{x}dt.
\end{equation*}%
Using (\ref{bo}), (\ref{W}), the convergence \ $\mathbf{1}_{\sigma }(\mathbf{%
x})\underset{\sigma \rightarrow 0}{\longrightarrow }1\ $\ \ for a.e.$\
\mathbf{x}\in \Omega ,$ we derive the inequality (\ref{A}).\bigskip

\subsection{\protect\vspace{1pt}Estimates independent of $\protect\theta .$
Limit transition on $\protect\theta $}

\bigskip \label{sec555}

In the previous sections we have constructed the solution $\{\omega ,h,%
\mathbf{v}=\mathrm{rot}\;h\}$ for the system (\ref{eq7}), (\ref{eq4sec2}),\
depending on ${R,\theta .}$

Let us assume that $R:=\frac{1}{\theta }.$ Therefore we can consider that $%
\left\{ \omega ,h, \mathbf{v} \right\} $ depends only on the
{ parameter} $\theta
.$ In the sequel instead of $\omega ,h,h_{1}$ (see the representation (\ref%
{Apresent})), $\mathbf{v}$ and the approximated data $\gamma \,\mathbf{,}%
g,A,a,$ we shall write $\omega _{\theta },$ $h_{\theta },h_{1,\theta },%
\mathbf{v_{\theta }}\ $ and $\gamma ^{\theta }\,\mathbf{,}g^{\theta
},A^{\theta },a^{\theta },$ respectively. 

\textit{In this subsection all
constants }$C$\textit{\ }\underline{\textit{do not depend on } $\theta $}.

\vspace{1pt}Let us formulate for Gronwall's type lemma. This lemma can be
proved by the standard \ method.\textbf{\ }

\begin{lemma}
\label{lem4.0} Let $D(t),B(t)\in L_{1}(0,T)$\ be given
{ non-negative } functions. Let $y(t)$ be a
non-negative function for $t\in \lbrack 0,T]$ and $y(t)=0$ for
$\forall t>T,$ satisfying
\begin{equation}
y(t)\leqslant y_{0}+\int_{0}^{t}\left[ D(s)\cdot u(s)+B(s)\right] ds
\label{gronwall}
\end{equation}%
with $u(t):=\frac{1}{\theta }\int_{t}^{t+\theta }y(s)\,ds.$ \ Then there
exists $\theta _{0}>0,$ such that for any fixed $\theta \in (0,$ $\theta _{0})$
we have
\begin{equation*}
y(t)\leqslant 2 \exp \left( \int_{0}^{t}D(s)ds\right) \left[
y_{0}+\int_{0}^{t}B(r)\cdot \exp \left( -\int_{0}^{r}D(s)ds\right) dr\right]
,\quad \quad \forall t\in \lbrack 0,T].
\end{equation*}
\end{lemma}

\bigskip

Combining Theorem \ref{lem6sec4234} and\ Lemma \ref{lem4.0}, we show the
following result.

\begin{lemma}
\label{teo4sec2 copy(1)} The pair $\{\omega _{\theta },h_{\theta }\}$
satisfies the estimates
\begin{align}
||\omega _{\theta }||_{L_{\infty }(0,T;\;L_{p}(\Omega ))}& \leqslant C,
\label{t1} \\
||h_{\theta }||_{L_{\infty }(0,T;\;W_{p}^{2}(\Omega ))}& \leqslant C,
\label{t2} \\
\left\Vert  \nabla _{\mathbf{x}%
}h_{1,\theta } (\cdot , t+\Delta )-  \nabla _{\mathbf{x}%
}h_{1,\theta }  (\cdot , t )  \right\Vert _{ L_{p}(\Omega ))}& \leqslant C \, \Delta   \label{t3}
\end{align}
for a.e. $t \in (0,T)$ and any  $ | \Delta | < \min \{t, T-t \} .$  
\end{lemma}

\noindent

\textbf{Proof. }  From (\ref{boundary}), (\ref{Bpresent}) and (\ref{ell-5}), (\ref{ell-8}) of Lemma \ref{teo2sec2},
 we have%
\begin{eqnarray}
|\omega _{\Gamma }(\mathbf{v}_{\theta })| &\leqslant &||\mathbf{v}_{\theta
}||_{L_{\infty }(0,T;C(\overline{\Omega }))}\cdot |\gamma ^{\theta }|\
+|g^{\theta }|\quad \quad \mbox{  on   }\quad \Gamma _{T}^{-},  \label{vvv}
\\
||\mathbf{v}_{\theta }(\cdot,t)||_{C(\overline{\Omega })}^{p}
&\leqslant &C_p\left\{ ||\left\langle \omega _{\theta }\right\rangle (\cdot
,t)||_{L_{p}(\Omega )}^{p}+{}||A^{\theta }(\cdot
,t)||_{W_{p}^{2}(\Gamma )}^{p}\right\} ,  \label{VV}
\end{eqnarray}%
for a.e. $t\in \lbrack 0,T]$. Here and below we denote by $C_p$ constants which  depend on $p$, being indepedent of $\theta$.
With the help of (\ref{eq122sec2}), we obtain
\begin{equation}
\int_{0}^{t}\int_{\Gamma ^{-}}a^{\theta }|\omega _{\Gamma }(\mathbf{v}%
_{\theta })|^{p}d\mathbf{x}d\tau \leqslant \int_{0}^{t}\left\{ D(s)\cdot
\frac{1}{\theta }\int_{s}^{s+\theta }||\omega _{\theta }(\cdot ,\tau
)||_{L_{p}(\Omega )}^{p}d\tau +B(s)\right\} \,ds,  \label{ss}
\end{equation}%
\begin{eqnarray*}
D(t) &:&=C_{p}||a||_{L_{\infty }(0,T;C(\Gamma ^{-}))}\cdot \int_{\Gamma
^{-}}\,|\gamma |^{p}d\mathbf{x}, \\
B(t) &:&=C_{p}||a||_{L_{\infty }(0,T;C(\Gamma ^{-}))}\cdot \left(
||A||_{L_{\infty }(0,T;W_{p}^{2}(\Gamma ))}^{p}\cdot \int_{\Gamma
^{-}}\,|\gamma |^{p}d\mathbf{x}+\int_{\Gamma ^{-}}\,|g|^{p}\,d\mathbf{x}%
\right) .
\end{eqnarray*}
 Taking into account (\ref{A}), (\ref{ss}) and applying Lemma \ref{lem4.0} to
the function $y(t):=||\omega _{\theta }(\cdot ,t)||_{L_{p}(\Omega )}^{p},$
we derive the estimate
\begin{equation}
||\omega _{\theta }(\cdot ,t)||_{L_{p}(\Omega )}^{p}\leqslant 2%
\exp \left( \int_{0}^{t}D(s)ds\right) \times \left[ ||\omega
_{0}||_{L_{p}(\Omega )}^{p}+\int_{0}^{t}B(s)ds\right] ,  \label{AAA}
\end{equation}%
for a.e. $t\in \lbrack 0,T]$ and $\forall \theta \in (0,\min \left\{ \theta
_{0},\sigma _{0},\frac{T}{4}\right\} ),$ that gives \eqref{t1}.

According to the representation (\ref{Apresent}), Lemma \ref{teo2sec2} and %
\eqref{t1}, we have%
\begin{equation*}
||h_{\theta }(\cdot ,t)||_{W_{p}^{2}(\Omega )}\leqslant C\left(
||\left\langle \omega _{\theta }\right\rangle (\cdot ,t)||_{L_{p}(\Omega
)}+||a^{\theta }(\cdot ,t)||_{W_{p}^{1}(\Gamma )}\right) \quad \quad
\mbox{  for  a.e.
}\quad t\in \left[ 0,T\right] ,
\end{equation*}%
hence, by (\ref{AregA}), \eqref{eq122sec2}, \eqref{eq122sec2-2} and \eqref{t1}, we obtain (%
\ref{t2}).

Now let us derive the estimate (\ref{t3}) for the derivative of $\nabla _{%
\mathbf{x}}h_{_{1},\theta }$ with respect to the time $t.$ Using the theory
of \cite{dip}, the equation (\ref{eq7}) implies the equality%
\begin{equation}
\int_{\Omega _{T}} [\omega _{\theta }]_{\frac 1 \theta } (\psi _{t}+%
\mathbf{v}_{\theta }\,\nabla \psi )\;d\mathbf{y}dt=0  \label{R-11}
\end{equation}%
for any function $\psi \in C_{0}^{1,1}(\Omega _{T}).$ It is obvious that in
(\ref{R-11}) we can put 
$$
\psi (\mathbf{x},t)=K_{1}(\mathbf{x},%
\mathbf{y})  \cdot  \int_t^{t+\Delta } \varphi (s)\, ds 
$$
with {\it arbitrary}  $\varphi \in C(0,T)$, having $supp (\varphi ) \subset (0,T)$, and deduce the
following equality
\begin{equation*}
h_{_{1},\theta }  (\mathbf{x}%
,t+\Delta )-h_{_{1},\theta } (\mathbf{x}%
,t)= \int_{\Omega }  \int_{t-\Delta}^t \left( [\omega _{\theta } (%
\mathbf{y}, s ) ]_{\frac 1 \theta } 
\ \mathbf{v}_{\theta }(\mathbf{y}, s)  \right) \, ds\cdot \nabla _{y}K_{1}(%
\mathbf{x},\mathbf{y})\;d\mathbf{y} 
\end{equation*}
for a.e. $t \in (0,T)$ and any  $ | \Delta | < \min \{t, T-t \} .$ 
Calculating the first derivative
with respect to $\mathbf{x}$, applying the properties (\ref{KK}) of the
kernel $K_{1}(\mathbf{x},\mathbf{y})$ and accounting (\ref{t1}), (\ref{t2}),  we derive 
\begin{equation*}
\left\Vert 
\nabla _{\mathbf{x}} h_{_{1},\theta }  (\cdot,t+\Delta )-\nabla _{\mathbf{x} } h_{_{1},\theta } (\cdot ,t)
\right\Vert _{L_{p}(\Omega )}\leqslant
C\left\Vert \mathbf{v}_{\theta } \right\Vert _{L_\infty (0,T;\, C(\overline{\Omega }%
)} \,\, \int_{t-\Delta}^t \left\Vert [\omega _{\theta } (%
\cdot , s ) ]_{\frac 1 \theta }   \right\Vert _{L_{p}(\Omega )} \, ds \, \leqslant C \Delta .   
\end{equation*}%
  $\blacksquare $\bigskip

From (\ref{t1})-(\ref{t3}), the representation (\ref{Apresent}) and the
approximation convergence \ (\ref{eq122sec2-2}), we conclude that there
exists a subsequence of $\{\omega _{\theta },h_{\theta },\mathbf{v}_{\theta
}\},$ such that
\begin{align}
h_{\theta }& \rightharpoonup h\quad \quad \mbox{ weakly}-\ast \mbox{ in }%
L_{\infty }(0,T;\,W_{p}^{2}(\Omega )),  \notag \\
\omega _{\theta }& \rightharpoonup \omega \quad \quad \mbox{ weakly}-\ast
\mbox{ in
}L_{\infty }(0,T;L_{p}(\Omega )),  \label{transition} \\
\mathbf{v}_{\theta }& \rightarrow \mathbf{v:=}\mathrm{rot}\;h\quad \quad
\mbox{   in
}L_{\infty }(0,T;L_{p}(\Omega )).  \notag
\end{align}%
Hence the limit transition on $\theta \rightarrow 0$ \ in the equations (\ref%
{eq7}) and (\ref{eq4sec2}), written for $\{\omega _{\theta },h_{\theta },%
\mathbf{v}_{\theta }\}$\vspace{1pt}, implies that the limit triple $\{\omega
,h,\mathbf{v}\}$ satisfies (\ref{eq7})-(\ref{eq10}).

\bigskip

Above we considered the case when $p\in (2,\infty ).$

Now we attend to the case $p=\infty .$ We assume that a set of approximate
data $\gamma ^{\theta },$ $g^{\theta },$ $A^{\theta },$\ $\omega
_{0}^{\theta },$ $a^{\theta }$\ satisfies the conditions (\ref{reg0})-(\ref%
{eq122sec2-2}), written for any $q\in (2,\infty )$ (instead of $p).$\
Additionally we suppose that (\ref{eq122sec2}) is true for $q=\infty .$ By
the above construction, there exists an approximate solution $\{\omega
_{\theta },h_{\theta },\mathbf{v}_{\theta }=\mathrm{rot}\;h_{\theta }\},$
satisfying all estimates of Lemma \ref{teo4sec2 copy(1)}. These estimates
are valid for any $q\in (2,\infty )$ (instead of given $p).$ By Theorem \ref%
{lem6sec4234} we have
\begin{equation}
||\omega _{\theta }(\mathbf{x},t_{0})||_{L_{q}(\Omega )}\leqslant ||\omega
_{0}^{\theta }||_{L_{q}(\Omega )}+\left( \int_{0}^{t_{0}}\int_{\Gamma
^{-}}a^{\theta }\;|\omega _{\Gamma }(\mathbf{v}_{\theta })|^{q}\;d\mathbf{x}%
dt\right) ^{1/q}  \label{inf2}
\end{equation}%
for all $\;t_{0}\in \lbrack 0,T]$ and $\forall q\in (2,\infty ).$\ Next we
use the estimate (\ref{VV}), written for some fixed $p=\overline{q}\in
(2,\infty ),$%
\begin{equation}
||\mathbf{v}_{\theta }||_{L_{\infty }(0,T;C(\overline{\Omega }))}\leqslant
C\left\{ ||\left\langle \omega _{\theta }\right\rangle (\cdot
,t)||_{L_{\infty }(0,T;L_{\overline{q}}(\Omega ))}+{}||A^{\theta }(\cdot
,t)||_{L_{\infty }(0,T;W_{\overline{q}}^{2}(\Gamma ))}\right\} \leqslant C(%
\overline{q}),  \label{inf1}
\end{equation}%
where $C(\overline{q})$ depends only on $\overline{q}.$ Then, applying (\ref%
{eq00sec1}), (\ref{eq122sec2}), (\ref{vvv}) in the inequality (\ref{inf2}),
we derive%
\begin{equation*}
||\omega _{\theta }(\mathbf{x},t_{0})||_{L_{q}(\Omega )}\leqslant ||\omega
_{0}||_{L_{q}(\Omega )}+C(\overline{q}).
\end{equation*}%
Passing to the limit on $q\rightarrow \infty ,\ $we$\ $come to the estimate%
\begin{equation}
||\omega _{\theta }(\mathbf{x},t_{0})||_{L_{\infty }(\Omega )}\leqslant
||\omega _{0}||_{L_{\infty }(\Omega )}+C(\overline{q})=:R_{0}.  \label{max}
\end{equation}%
Choosing $R=\frac{1}{\theta }>R_{0},$ we can take off the subscript $R$ in (%
\ref{R}) and in remaining formulas. The derived estimates permit to pass to
the limit on $\theta \rightarrow 0$\ in (\ref{eq7}) and (\ref{eq4sec2}),
written for $\{\omega _{\theta },h_{\theta },\mathbf{v}_{\theta }\},$\ and
prove the existence of the weak solution of the problem \eqref{eq1001}-%
\eqref{eq1006}, \eqref{eq1.39}-\eqref{eqC2}, \eqref{eqq1.2} in the case $p=\infty$.

This concludes the proof of Theorem \ref{teo2sec1}.

\bigskip

\section*{Acknowledgement}

N.V. Chemetov thanks for support from FCT and FEDER through the
Project POCTI / ISFL / 209 of Centro de Matem\'{a}tica e
Aplica\c{c}\~{o}es Fundamentais da Universidade de Lisboa (CMAF /
UL) and the Project POCTI / MAT / 45700 / 2002. { The
work of S.N. Antontsev partially  supported  by the Portuguese
research project POCI/MAT/61576/2004/FCT/MCES.}

\end{document}